\documentclass[a4paper]{article}
\usepackage{amsmath,amsthm,amssymb,amscd}

\newcommand{\corresp}[3]{\raisebox{-1ex}[0ex][1.3ex]{\scriptsize $#2$}\hspace{0.3ex}
  \framebox{$#1$}\hspace{0.3ex}
  \raisebox{-1ex}[0ex][1.3ex]{\scriptsize $#3$}}
\newcommand{\bicorresp}[4]{\raisebox{-1ex}[0ex][1.3ex]{\scriptsize $#2$}\hspace{0.3ex}
  \overset{#4}{\framebox{$#1$}}\hspace{0.3ex} \raisebox{-1ex}[0ex][1.3ex]{\scriptsize $#3$}}

\newcommand{\pil}{\pi_{\ell}}
\newcommand{\pitill}{\widetilde{\pi}_{\ell}}
\newcommand{\pir}{\pi_{r}}
\newcommand{\Ind}{\operatorname{Ind}}

\newcommand{\Ctil}{\widetilde{C}}
\newcommand{\tetil}{\widetilde{\theta}}

\newcommand{\Vh}{\hat{V}}
\newcommand{\Yh}{\hat{Y}}
\newcommand{\M}{\mathcal{M}}
\newcommand{\B}{\operatorname{B}}
\newcommand{\K}{\mathcal{K}}
\newcommand{\E}{\mathcal{E}}
\newcommand{\Etil}{{\widetilde{\E}}}
\newcommand{\F}{\mathcal{F}}
\newcommand{\Ftil}{{\widetilde{\F}}}
\newcommand{\I}{\mathcal{I}}
\newcommand{\J}{\mathcal{J}}
\newcommand{\G}{\mathcal{G}}
\newcommand{\la}{\langle}
\newcommand{\ra}{\rangle}
\newcommand{\te}{\theta}

\newcommand{\ltimesred}{\ {_{\fontshape{n}\selectfont
  r}}\hspace{-0.8ex}\ltimes}
\newcommand{\ltimesredsmall}{\, {_{\fontshape{n}\selectfont
  r}}\hspace{-0.4ex}\ltimes}
\newcommand{\ltimesfull}{\ {_{\fontshape{n}\selectfont
  f}}\hspace{-0.8ex}\ltimes}
\newcommand{\ltimesfullsmall}{\, {_{\fontshape{n}\selectfont
  f}}\hspace{-0.4ex}\ltimes}
\newcommand{\rtimesred}{\rtimes_{\fontshape{n}\selectfont r}}
\newcommand{\rtimesfull}{\rtimes_{\fontshape{n}\selectfont f}}

\newcommand{\Morita}{\underset{\text{\rm Morita}}{\sim}}

\newcommand{\Efull}{\mathcal{E}_{\fontshape{n}\selectfont
  f}}
\newcommand{\Eprimefull}{\mathcal{E}'_{\fontshape{n}\selectfont
  f}}

\newcommand{\Ahu}{{\hat{A}}^{\text{\rm u}}}
\newcommand{\Au}{A^{\text{\rm u}}}
\newcommand{\Xu}{X_{\text{\rm u}}}
\newcommand{\piu}{\pi_{\text{\rm u}}}
\newcommand{\deu}{\Delta^{\text{\rm u}}}
\newcommand{\dehu}{{\hat{\de}}^{\text{\rm u}}}
\newcommand{\Wu}{W^{\text{\rm u}}}
\newcommand{\etah}{\hat{\eta}}

\newcommand{\rot}[1]{\underset{#1}{\ot}}
\newcommand{\rotalg}[1]{\overset{\operatorname{alg}}{\underset{#1}{\ot}}}

\newcommand{\Ah}{\hat{A}}
\newcommand{\Rh}{\hat{R}}

\newcommand{\cst}{C$^*$}
\newcommand{\wst}{W$^*$}
\newcommand{\pih}{\hat{\pi}}
\newcommand{\vfi}{\varphi}
\newcommand{\Mfi}{\mathcal{M}_\varphi}
\newcommand{\Mpsi}{\mathcal{M}_\psi}
\newcommand{\Nfi}{\mathcal{N}_\varphi}
\newcommand{\cC}{\mathcal{C}}

\newcommand{\vfih}{\hat{\varphi}}

\newcommand{\Ahop}{\hat{A}\hspace{-.1ex}\raisebox{0.9ex}[0pt][0pt]{\scriptsize\fontshape{n}\selectfont op}}
\newcommand{\Ahopu}{\Ahu\hspace{-.2ex}\raisebox{0.9ex}[0pt][0pt]{\scriptsize\fontshape{n}\selectfont ,op}}
\newcommand{\Aop}{A\hspace{-.2ex}\raisebox{0.9ex}[0pt][0pt]{\scriptsize\fontshape{n}\selectfont op}}
\newcommand{\Aopu}{\Au\hspace{-.2ex}\raisebox{0.9ex}[0pt][0pt]{\scriptsize\fontshape{n}\selectfont ,op}}

\newcommand{\etatil}{\widetilde{\eta}}
\newcommand{\gatil}{\widetilde{\gamma}}

\newcommand{\alh}{\hat{\alpha}}
\newcommand{\Mh}{\hat{M}}
\newcommand{\Si}{\Sigma}

\newcommand{\dpr}{^{\prime\prime}}

\newcommand{\Wh}{\hat{W}}

\newcommand{\C}{\mathbb{C}}

\newcommand{\Q}{\mathbb{Q}}

\newcommand{\deh}{\hat{\Delta}}

\newcommand{\cJ}{{\cal J}}

\newcommand{\Jh}{\hat{J}}

\newcommand{\cI}{{\cal I}}
\newcommand{\cL}{{\cal L}}
\newcommand{\cF}{{\cal F}}

\newcommand{\cN}{{\cal N}}
\newcommand{\cM}{{\cal M}}

\newcommand{\ot}{\otimes}

\newcommand{\om}{\omega}
\newcommand{\io}{\iota}
\newcommand{\al}{\alpha}
\newcommand{\be}{\beta}
\newcommand{\ga}{\gamma}

\newcommand{\de}{\Delta}
\newcommand{\si}{\sigma}
\newcommand{\cV}{{\cal V}}
\newcommand{\cW}{{\cal W}}

\newcommand{\deop}{\de \hspace{-.3ex}\raisebox{0.9ex}[0pt][0pt]{\scriptsize\fontshape{n}\selectfont op}}
\newcommand{\deopo}{\de_1
\hspace{-1ex}\raisebox{1ex}[0pt][0pt]{\scriptsize\fontshape{n}\selectfont
op}}

\newcommand{\dehop}{\deh \hspace{-.3ex}\raisebox{0.9ex}[0pt][0pt]{\scriptsize\fontshape{n}\selectfont op}}
\newcommand{\dehopo}{\deh_1
\hspace{-1ex}\raisebox{1ex}[0pt][0pt]{\scriptsize\fontshape{n}\selectfont
op}}

\newcommand{\recht}{\rightarrow}

\numberwithin{equation}{section}

{\theoremstyle{definition}\newtheorem{definition}{Definition}[section]
\newtheorem{notation}[definition]{Notation}

\newtheorem{remark}[definition]{Remark}
}
\newtheorem{proposition}[definition]{Proposition}
\newtheorem{lemma}[definition]{Lemma}
\newtheorem{theorem}[definition]{Theorem}
\newtheorem{corollary}[definition]{Corollary}

\begin{document}
\begin{center}
{\LARGE\bf A new approach to induction and imprimitivity results}

\bigskip

{\sc by Stefaan
  Vaes}
\end{center}

{\footnotesize
Institut de Math{\'e}matiques de Jussieu; Alg{\`e}bres d'Op{\'e}rateurs et
Repr{\'e}sentations; 175, rue du Chevaleret; F--75013 Paris (France) \\
Department of Mathematics; K.U.Leuven; Celestijnenlaan 200B; B-3001 Leuven (Belgium) \\
e-mail: vaes@math.jussieu.fr}

\bigskip

\begin{abstract}\noindent
In the framework of locally compact quantum groups, we provide an induction procedure for unitary corepresentations as well as coactions on
\cst-algebras. We prove imprimitivity theorems that unify the existing theorems for actions and coactions of groups. We essentially use von Neumann
algebraic techniques.
\renewcommand{\thefootnote}{}
\end{abstract}

\section{Introduction}

The theory of induced representations of locally compact (l.c.) groups
was introduced by Mackey \cite{Mackey}, who discovered the
\emph{imprimitivity theorem}, characterizing induced representations
through the presence of a covariant representation of a homogeneous
space. Rieffel \cite{Rieffel} provided a modern approach using the
language of Hilbert \cst-modules.

After the work of Rieffel, several induction procedures and
imprimitivity results have been obtained, both for actions and for
coactions of groups. We shall briefly review them below.
The purpose of this paper is to develop such an induction and
imprimitivity machinery in the setting of \emph{locally compact
  quantum groups} \cite{KV1,KV2}. In this way, we provide a unified
approach to several results on actions and coactions of groups
on \cst-algebras. At
the same time, our proofs in the setting of l.c.\ quantum groups are
simpler than the classical proofs dealing with coactions of groups.

Another motivation comes from quantum group theory. Meyer and Nest
have undertaken a reformulation of the Baum-Connes conjecture
\cite{MN} in which induction and restriction play a crucial role. The
development of induction and imprimitivity in the current paper should
play an equally important role in the formulation of a Baum-Connes conjecture
for quantum groups.

A final aspect of the paper is the technique that is used to prove the imprimitivity results: we shall obtain Morita equivalences between
\cst-algebras by using von Neumann algebra techniques and the language of correspondences \cite{CJ,Popa}. This is the main reason why our proofs are
simpler than the classical proofs dealing with coactions. Already for group duals, but certainly for l.c.\ quantum groups, the von Neumann algebra
picture of the quantum group (looking at $L^\infty$ rather then $C_0$) is much more user-friendly.

A first approach to induction of unitary corepresentations of l.c.\ quantum groups has been developed by Kustermans \cite{JK2}, but without dealing
with imprimitivity results or coactions on \cst-algebras. Of course, our induction procedure is unitarily equivalent to his. Our approach is simpler
and makes it more easy to prove properties (e.g.\ induction in stages).

In the remainder of the introduction, we shall review several
imprimitivity results that were obtained for actions and coactions of
groups. We shall explain how they are generalized in the setting of
l.c.\ quantum groups.

Rieffel has given a modern \cst-algebraic formulation of Mackey's result using Hilbert \cst-modules \cite{Rieffel}. The neatest form of Mackey's
result is given by the Morita equivalence of \cst-algebras
\begin{equation} \label{eq.mackey}
G \ltimesfull C_0(G/G_1) \Morita C^*(G_1) \; ,
\end{equation}
whenever $G_1$ is a
closed subgroup of a l.c.\ group $G$.

Green \cite{Green} has generalized Mackey's induction of unitary representations to
\cst-dynamical systems. Suppose that $G_1$ is a closed subgroup of a
l.c.\ group $G$. Suppose that $G_1$ acts continuously on a
\cst-algebra $B$. Then, Green constructs an induced \cst-algebra $\Ind
B$ with a continuous action of $G$ such that we obtain a Morita equivalence
\begin{equation}\label{eq.green}
G \ltimesfull \Ind B \Morita
G_1 \ltimesfull B \; .
\end{equation}
If $B=\C$, we get $\Ind
B = C_0(G/G_1)$ and we find back Mackey's imprimitivity theorem \eqref{eq.mackey}.

If a l.c.\ group $G$ acts continuously on a \cst-algebra $B$ and if $G_1$ is a
closed subgroup of $G$, we can first restrict the action of $G$ to an
action of $G_1$ on $B$ and then induce this restricted action to an
action of $G$. The resulting \cst-algebra is $C_0(G/G_1) \ot B$ with the
diagonal action of $G$. If we now suppose that $G_1$ is normal in $G$,
Green's imprimitivity theorem can be restated as the Morita
equivalence
\begin{equation}\label{eq.mackeygreen}
G_1 \ltimesfull B \Morita \widehat{G/G_1} \ltimes (G \ltimesfull B) \; .
\end{equation}
The second crossed
product is the crossed product by the restriction of the dual coaction
to $G/G_1$.

Dually, Mansfield proved a coaction version of the Morita equivalence
\eqref{eq.mackeygreen}: if $G_1$ is a closed normal subgroup of a
l.c.\ group $G$ and if $B$ is a \cst-algebra with a \emph{reduced
coaction} of $G$, we have the following Morita equivalence for
\emph{reduced crossed products}:
\begin{equation}\label{eq.mansfield}
\widehat{G/G_1} \ltimes B \Morita G_1 \ltimesred (\widehat{G} \ltimes
B) \; .
\end{equation}
In fact, Mansfield had to impose the amenability of $G_1$ and Kaliszewski
\& Quigg \cite{KQ1} showed the result for general normal subgroups
$G_1$. The terminology \emph{reduced coaction} is not the standard
one: in the literature one uses \emph{normal coaction} (see Definition \ref{def.reducedcoaction}).

With representation theory in mind, one wants of course imprimitivity
results between full crossed products. Kaliszewski and Quigg
have shown recently that \eqref{eq.mansfield} holds for full crossed
products and \emph{maximal coactions} of $G$ (as introduced by
Echterhoff, Kaliszewski and Quigg \cite{EKQ}). A maximal coaction is a coaction which is Morita equivalent to a dual
coaction on a full crossed product (and a reduced coaction is a
coaction which is Morita equivalent to a dual coaction on a reduced
crossed product), see Definition \ref{def.maximalcoaction}. The same imprimitivity result
had been proved before for dual coactions by Echterhoff, Kaliszewski and Raeburn \cite{EKR}.

In Section \ref{sec.qhs} we define the \emph{quantum homogeneous
  space}, given a l.c.\ quantum group and a closed quantum subgroup. We
  prove a quantum version of the Mackey imprimitivity theorem
  \eqref{eq.mackey}. In Section \ref{sec.indcoac}, we study dynamical
  systems. Given a coaction of a closed quantum subgroup $(A_1,\de_1)$
  of a l.c.\ quantum group $(A,\de)$ on a
  \cst-algebra $B$, we construct an induced \cst-algebra $\Ind B$ with
  a coaction of $(A,\de)$ such that a quantum version of
  \eqref{eq.green} holds. Observe that already for coactions of groups
  such a construction of induced coactions was not known up to now.

In Section \ref{sec.indres} we describe what happens if we first
restrict and then induce a coaction. Instead of a tensor product, we
obtain a twisted product of the original \cst-algebra and the quantum
homogeneous space with some kind of diagonal coaction. This is used in
Section \ref{sec.mansfield} to obtain a quantum version of
\eqref{eq.mansfield}, both for reduced crossed products (and reduced
coactions) and for full crossed products (and maximal coactions). In
fact we define crossed products by homogeneous spaces and hence, we do
not have to assume normality of the quantum subgroup.

In \cite{EKQR}, Echterhoff, Kaliszewski, Quigg and Raeburn discuss
naturality of the imprimitivity theorems for actions and coactions of
groups. In our general approach, we also get covariant Morita equivalences and naturality.

As stated above, the technique used in this paper is von Neumann
algebraic in nature. After a section of preliminaries on locally
compact quantum groups, closed quantum subgroups and crossed products, we provide a
von Neumann algebraic approach to representation theory for quantum
groups in Section \ref{sec.fourdiff}. This is used to give an easy
approach to induction of representations in Section
\ref{sec.induction}, simplifying the original approach by Kustermans
\cite{JK2}. In Section \ref{sec.firstimprim} we prove a preliminary
imprimitivity theorem, which is the crucial ingredient in the next
sections that we already discussed above.

\begin{notation}
The most over-used symbol of this paper is $\ot$. It shall be used to
denote tensor products of Hilbert spaces and von Neumann algebras, as
well as \emph{minimal} tensor products of \cst-algebras.

The multiplier algebra of a \cst-algebra $A$ is denoted by $\M(A)$.

When $X$ is a subset of a Banach space, we denote by $[X]$ the closed linear span of $X$.
\end{notation}

\section{Preliminaries}

\subsection{Locally compact quantum groups}
We use \cite{KV1,KV2} as references for the theory of locally compact
(l.c.) quantum groups. Since the definition of a l.c.\ quantum group
in \cite{KV1,KV2} is based on the existence of invariant weights (Haar
measures), we introduce the following weight theoretic notation.

Let $\vfi$ be a normal, semi-finite, faithful (n.s.f.) weight on a von Neumann algebra $M$. Then, we write
$$
\cM_\vfi^+ = \{ x \in M^+ \mid \vfi(x) < \infty \} \quad\text{and}
\quad \cN_\vfi = \{ x \in M \mid \vfi(x^* x) < \infty \} \; .
$$

\begin{definition}
A pair $(M,\de)$ is called a (von Neumann algebraic) l.c.\ quantum group when
\begin{itemize}
\item $M$ is a von Neumann algebra and $\de : M \recht M \ot M$ is
a normal and unital $*$-homomorphism satisfying the coassociativity relation : $(\de \ot \io)\de = (\io \ot \de)\de$;
\item there exist n.s.f.\ weights $\varphi$ and $\psi$ on $M$ such that
\begin{itemize}
\item $\varphi$ is left invariant in the sense that $\varphi \bigl( (\om \ot
\io)\de(x) \bigr) = \varphi(x) \om(1)$ for all $x \in \Mfi^+$ and $\om \in M_*^+$,
\item $\psi$ is right invariant in the sense that $\psi \bigl( (\io \ot
\om)\de(x) \bigr) = \psi(x) \om(1)$ for all $x \in \Mpsi^+$ and $\om \in M_*^+$.
\end{itemize}
\end{itemize}
\end{definition}
Fix such a l.c.\ quantum group $(M,\de)$.

Represent $M$ in the GNS-construction of $\vfi$ with GNS-map $\Lambda : \Nfi \recht H$. We define a unitary $W$ on $H \ot H$ by
$$W^* (\Lambda(a) \ot \Lambda(b)) = (\Lambda \ot \Lambda)(\de(b)(a \ot 1)) \quad\text{for all}\; a,b \in \cN_{\vfi}
\; .$$ Here, $\Lambda \ot \Lambda$ denotes the canonical GNS-map for the tensor product weight $\varphi \ot \varphi$. One proves that $W$ satisfies
the pentagonal equation: $W_{12} W_{13} W_{23} = W_{23} W_{12}$, and we say that $W$ is a \emph{multiplicative unitary}. It is the \emph{left regular
corepresentation}. The von Neumann algebra $M$ is the strong closure of the algebra $\{ (\io \ot \om)(W) \mid \om \in \B(H)_* \}$ and $\de(x) = W^*
(1 \ot x) W$, for all $x \in M$. Next, the l.c.\ quantum group $(M,\de)$ has an antipode $S$, which is the unique $\si$-strong$^*$ closed linear map
from $M$ to $M$ satisfying $(\io \ot \om)(W) \in D(S)$ for all $\om \in \B(H)_*$ and $S(\io \ot \om)(W) = (\io \ot \om)(W^*)$ and such that the
elements $(\io \ot \om)(W)$ form a $\si$-strong$^*$ core for $S$. The antipode $S$ has a polar decomposition $S = R \tau_{-i/2}$ where $R$ is an
anti-automorphism of $M$ and $(\tau_t)$ is a strongly continuous one-parameter group of automorphisms of $M$. We call $R$ the \emph{unitary antipode}
and $(\tau_t)$ the \emph{scaling group} of $(M,\de)$. From \cite{KV1}, Proposition~5.26 we know that $\sigma (R \ot R) \de = \de R$. So $\varphi R$
is a right invariant weight on $(M,\de)$ and we take $\psi:= \varphi R$.

The \emph{dual l.c.\ quantum group} $(\Mh,\deh)$ is defined in \cite{KV1}, Section~8. Its von Neumann algebra $\Mh$ is the strong closure of the
algebra $\{(\om \ot \io)(W) \mid \om \in \B(H)_* \}$ and the comultiplication is given by $\deh(x) = \Sigma W (x \ot 1) W^* \Sigma$ for all $x \in
\Mh$. On $\Mh$ there exists a canonical left invariant weight $\vfih$ and the associated multiplicative unitary
is denoted by $\hat{W}$. From \cite{KV1}, Proposition~8.16, it follows that $\hat{W} = \Sigma W^* \Sigma$.

Since $(\Mh,\deh)$ is again a l.c.\ quantum group, we can introduce
the antipode $\hat{S}$, the unitary antipode $\hat{R}$ and the scaling group
$(\hat{\tau}_t)$ exactly as we did it for $(M,\de)$.

We shall denote the modular conjugations of the weights $\varphi$ and
$\hat\varphi$ by $J$ and $\Jh$ respectively. The operators $J$ and
$\Jh$ are anti-unitary involutions on the Hilbert space $H$. They
implement the unitary antipodes in the sense that
$$R(x) = \Jh x^* \Jh \quad\text{for all} \; x \in M
\qquad\text{and}\qquad \Rh(y) = J y^* J \quad\text{for all}\; y \in \Mh \; .$$
From modular theory, we also know that $M' = JMJ$ and $\Mh' = \Jh \Mh \Jh$.

We already discussed the left regular corepresentations $W$ and $\Wh$ of $(M,\de)$ and $(\Mh,\deh)$, respectively. These are multiplicative unitaries
and $\Wh = \Si W^* \Si$. We observe moreover that $W \in M \ot \Mh$ and $\Wh \in \Mh \ot M$. Since we also have right invariant weights on $(M,\de)$
and $(\Mh,\deh)$, we consider as well the right regular corepresentations $V \in \Mh' \ot M$ and $\hat{V} \in M' \ot \Mh$. These are also
multiplicative unitaries and satisfy
$$V = (\Jh \ot \Jh) \Wh (\Jh \ot \Jh) \; , \quad \Vh = (J \ot J) W (J
\ot J) \; .$$
We finally mention the formula $W^* = (\Jh \ot J)W(\Jh \ot J)$, which
is equivalent to saying that $(R \ot \Rh)(W) = W$.

Every l.c.\ quantum group has an \emph{associated \cst-algebra} of \lq continuous functions tending to zero at infinity\rq\ and we denote it by $A$
(resp.\ $\Ah$):
$$A := [(\io \ot \om)(W) \mid \om \in \B(H)_* ] \; , \quad \Ah = [(\om
\ot \io)(W) \mid \om \in \B(H)_*] \; .$$

It is clear that $A \subset M$. Also, the comultiplication $\de$
restricts to $A$ and yields a comultiplication $\de : A \recht \M(A
\ot A)$.

\begin{notation}
When we shall be dealing with coactions of l.c.\ quantum groups on
\cst-algebras, we will make use all the time of the \cst-algebraic
picture $(A,\de)$ of our l.c.\ quantum group. So, we shall speak about
\emph{the l.c.\ quantum group $(A,\de)$}.
\end{notation}

\subsection{Closed quantum subgroups}

We first discuss the notion of a \emph{morphism} between l.c.\ quantum
groups: see the work of Kustermans \cite{JK} for details. We explained
that every l.c.\ quantum groups admits a \cst-algebra $A$ and a dual
\cst-algebra $\hat{A}$. In the classical case of l.c.\ groups, this
comes down to the \cst-algebras $C_0(G)$ and $C^*_r(G)$. But there is
of course as well the universal \cst-algebra $C^*(G)$.

\begin{definition}
Let $(A,\de)$ be a l.c.\ quantum group. A \emph{unitary
  corepresentation} of $(A,\de)$ on a \cst-$B$-module $\E$ is a
  unitary $X \in \cL(A \ot \E)$ satisfying $$(\de \ot \io)(X) =
  X_{13}X_{23} \; .$$
\end{definition}

In \cite{JK} the universal dual $(\Ahu,\dehu)$ of $(A,\de)$ is
defined. This definition is analogous to the definition of the full group \cst-algebra $C^*(G)$
of a l.c.\ group $G$. This means that there exists a
\emph{universal corepresentation} $\cW \in \M(A \ot \Ahu)$ such that
the formula
$$(\io \ot \te)(\cW) = X$$
gives a bijective correspondence between non-degenerate
$^*$-homomorphisms $\te : \Ahu \recht \cL(\E)$ and unitary
corepresentations $X \in \cL(A \ot \E)$, whenever $\E$ is a
\cst-$B$-module. Moreover, the comultiplication $\dehu$ satisfies
$$(\io \ot \dehu)(\cW) = \cW_{13} \cW_{12} \; .$$

In exactly the same way, there is a universal version of $(A,\de)$,
which is denoted by $(\Au,\deu)$. There exists a universal
corepresentation $\hat{\cW} \in \M(\Au \ot \Ah)$ of $(\Ah,\deh)$ such
that the formula $(\te \ot \io)(\hat{\cW}) = X$ gives a
bijective correspondence between non-degenerate
$^*$-homomorphisms $\te : \Au \recht \cL(\E)$ and unitary
corepresentations $X \in \cL(\E \ot \Ah)$, whenever $\E$ is a
\cst-$B$-module.

\begin{definition}
A morphism $(M,\de) \overset{\pi}{\longrightarrow} (M_1,\de_1)$
between the l.c.\ quantum groups $(M,\de)$ and $(M_1,\de_1)$ is a
non-degenerate $^*$-homomorphism
$$\pi : \Au \recht \M(\Au_1) \quad\text{satisfying}\quad \deu_1 \pi =
(\pi \ot \pi)\deu \; .$$
\end{definition}

We slightly abuse notation by writing $(M,\de)
\overset{\pi}{\longrightarrow} (M_1,\de_1)$ and we should always keep
in mind that $\pi$ lives on the level of universal \cst-algebras.

Observe that in the classical situation, this corresponds to $\Au =
C_0(G)$, $\Au_1 = C_0(G_1)$ and $\pi(f) = f\circ \te$, where $\te :
G_1 \recht G$ is a continuous group homomorphism. Associated with
$\te$, we can then write as well $\pih : C^*(G_1) \recht \M(C^*(G))$
defined by $\pih(\lambda_p) = \lambda_{\te(p)}$. In the same way,
every morphism $(M,\de) \overset{\pi}{\longrightarrow} (M_1,\de_1)$
between l.c.\ quantum groups admits canonically a dual morphism
$(\Mh_1,\deh_1) \overset{\pih}{\longrightarrow} (\Mh,\deh)$.

\begin{definition}
We say that the morphism $(M,\de) \overset{\pi}{\longrightarrow}
(M_1,\de_1)$ identifies $(M_1,\de_1)$ as a \emph{closed quantum
  subgroup} of $(M,\de)$ if there exists a faithful, normal, unital
$^*$-homomorphism $\Mh_1 \recht \Mh$ which makes the following diagram
commute.
$$\begin{CD}
\Ahu_1 @>{\pih}>> \M(\Ahu) \\
@VVV @VVV \\
\Mh_1 @>>> \Mh
\end{CD}$$
In that case, we continue writing $\pih : \Mh_1 \recht \Mh$.
\end{definition}

One verifies that, in the classical setting where $\pi(f) = f \circ
\te$ for a continuous group homomorphism $\te : G_1 \recht G$, this
comes down to the fact that $\te$ identifies $G_1$ with a closed
subgroup of $G$. Indeed, the map $\pih(\lambda_p) = \lambda_{\te(p)}$
extends to a faithful, normal $^*$-homomorphism $\cL(G_1) \recht
\cL(G)$ if and only if $G_1$ is a closed subgroup of $G$.

\subsection{Crossed products and regularity}

\begin{definition} \label{def.defcoac}
A \emph{coaction} of a l.c.\ quantum group $(A,\de)$ on a \cst-algebra
$B$ is a non-degenerate \linebreak $^*$-homomorphism
$$\al : B \recht \M(A \ot B) \quad\text{satisfying}\quad (\io \ot
\al)\al = (\de \ot \io)\al \; .$$
We say that $\al$ is a \emph{continuous coaction} if
$$[\al(B)(A \ot 1)] = A \ot B \; .$$
\end{definition}

Coactions and their associated crossed products have been studied in
detail by Baaj and Skandalis \cite{BS}. We recall some basic
concepts. In this paper, we shall make use as well of coactions on
\cst-modules. This is discussed in detail in the appendix, following
another paper of Baaj and Skandalis \cite{KKS}.

Let $\al : B \recht \M(A \ot B)$ be a continuous coaction of a l.c.\
quantum group $(A,\de)$ on a \cst-algebra $B$.
Then,
$$[\al(B) (\Ah \ot 1)] \subset \M(\K(H) \ot B)$$
is a \cst-algebra which is called the \emph{reduced crossed product}
and denoted by $\Ah \ltimesred B$.

A pair $(X,\pi)$ consisting of a unitary corepresentation $X \in \M(A
\ot \K(K))$ of $(A,\de)$ on a Hilbert space $K$ and a non-degenerate
$^*$-homomorphism $\pi : B \recht \B(K)$ is called a \emph{covariant
  representation} of $\al$ if
$$(\io \ot \pi)\al(x) = X^* (1 \ot \pi(x)) X \quad\text{for all}\quad
x \in B \; .$$ There exists a unique (up to isomorphism) \cst-algebra denoted by $\Ahu \ltimesfull B$ and called the \emph{full crossed product},
equipped with a universal covariant representation
$$\Xu \in \M(A \ot (\Ahu \ltimesfull B)) \; , \quad \piu : B \recht
\M(\Ahu \ltimesfull B)$$ such that the formulas
$$X = (\io \ot \te)(\Xu) \quad\text{and}\quad \pi = \te \piu$$
yield a bijective correspondence between covariant representations $(X,\pi)$ of $\al$ and non-degenerate representations $\te$ of the \cst-algebra
$\Ahu \ltimesfull B$.

By definition it is clear that $\Ahu$ coincides with the full crossed product of $(A,\de)$ coacting on the trivial \cst-algebra $\C$. It is also
clear that there is a natural surjective $^*$-homomorphism $\Ahu \ltimesfull B \recht \Ah \ltimesred B$.

\begin{remark}
From time to time we shall use as well right coactions $\al : B \recht
\M(B \ot A)$, satisfying $(\al \ot \io)\al = (\io \ot \de)\al$. The
reduced crossed product is then given by
$$B \rtimesred \Ahop = [\al(B) (1 \ot \Jh \Ah \Jh)] \; .$$
The reason why the opposite algebra $\Jh \Ah \Jh$ appears is natural: a right coaction of $(A,\de)$ corresponds to a left coaction of $(A,\deop)$.
The dual of the opposite quantum group $(A,\deop)$ is $\Jh \Ah \Jh$. In the same way, one defines $B \rtimesfull \Ahopu$.
\end{remark}

Both the reduced and the full crossed products admit a \emph{dual
  coaction} of $(\Ah,\dehop)$, which leaves invariant $B$ and acts as
  the comultiplication $\dehop$ on $\Ah$. So, it is a natural question
  to consider what happens with the second crossed products
$$\Aop \ltimesred (\Ah \ltimesred B) \quad\text{and}\quad \Aopu
\ltimesfull (\Ahu \ltimesfull B) \; .$$

For abelian l.c.\ groups, it is well known that $\hat{G} \ltimes G \ltimes B \cong \K(L^2(G)) \ot B$. This result need no longer be true for l.c.\
quantum groups. Indeed, taking $B = \C$, it might very well be the case that $\Aop \ltimesred \Ah \not\cong \K(H)$. The following definition, due to
Baaj and Skandalis \cite{BS}, describes the quantum groups for which the classical biduality result holds.

For any multiplicative unitary $W$ on a Hilbert space $H$, we introduce (following
\cite{BS}) the algebra $\cC(W)$ by the formula
$$\cC(W) := \{ (\io \ot \om)(\Si W) \mid \om \in \B(H)_* \} \; .$$

\begin{definition}
A l.c.\ quantum group $(A,\de)$ is said to be \emph{regular} if
$[\cC(W)] = \K(H)$, where $W$ denotes the left regular representation.
\end{definition}

It follows from Proposition 2.6 in \cite{BSV} that a l.c.\ quantum
group is regular if and only if the reduced crossed product of $A$ and
$\Ah$ is isomorphic with $\K(H)$.

\begin{remark}
The notion of \emph{continuous coaction} is somehow problematic for
non-regular quantum groups. Definition \ref{def.defcoac} makes sense, but is
not the only natural definition in the non-regular case. See \cite{BSV} for a
detailed discussion.
\end{remark}

\begin{definition} \label{def.reducedcoaction}
A continuous coaction $\al : B \recht \M(A \ot B)$ of $(A,\de)$ on $B$
is said to be \emph{reduced} if $\al$ is a faithful $^*$-homomorphism.
\end{definition}
Whenever $\al$ is a reduced continuous coaction of a regular l.c.\
quantum group $(A,\de)$ on a \cst-algebra $B$, we have that
$$\Aop \ltimesred (\Ah \ltimesred B) \cong \K(H) \ot B \; .$$
In fact, one obtains a natural covariant Morita equivalence $\Aop \ltimesred (\Ah \ltimesred B) \Morita B$ with respect to the bidual coaction on the
double crossed product and the coaction $\al$ on $B$.

All dual coactions on reduced crossed products are reduced coactions and the
biduality theorem shows that, in fact, a continuous coaction is
reduced if and only if it is Morita equivalent to a dual coaction.

\begin{definition} \label{def.strongreg}
A l.c.\ quantum group $(A,\de)$ is said to be \emph{strongly regular} if $\Ahu \ltimesfull A \cong \K(H)$.
\end{definition}

Since we always have the surjective $^*$-homomorphism $\pi: \Ahu \ltimesfull A \recht \Ah \ltimesred A$, it follows that a l.c.\ quantum group is
strongly regular if and only if it is regular and the $^*$-homomorphism $\pi$ is faithful.

\begin{remark}
Not every l.c.\ quantum group is regular. Non-regular examples are
given by the quantum groups $E_\mu(2)$, $ax+b$ or certain bicrossed
products (see \cite{BSV} for a detailed discussion of the latter
case). Examples of regular quantum groups include all the compact or
discrete quantum groups, all Kac algebras and a wide class of
bicrossed products. Also the analytic versions of the algebraic
quantum groups \cite{KVD} are regular.

All the above mentioned examples of regular quantum groups are in fact strongly regular. It is not known whether there exist regular quantum groups
which are not strongly regular.
\end{remark}

We choose to define coactions as $^*$-homomorphisms $\al : B \recht
\M(A \ot B)$, where $A$ is the \emph{reduced} \cst-algebra of the
l.c.\ quantum group. Of course, we can define
a continuous coaction of the universal quantum group $(\Au,\deu)$ on
the \cst-algebra $B$ as a non-degenerate $^*$-homomorphism $\al : B
\recht \M(\Au \ot B)$ satisfying $(\io \ot \al)\al = (\deu \ot
\io)\al$ and $[\al(B)(\Au \ot 1)] = \Au \ot B$. Observe that it
follows automatically that $(\varepsilon \ot \io)\al(x) = x$ for all
$x \in B$, where $\varepsilon : \Au \recht \C$ denotes the co-unit of
$(\Au,\deu)$. If we denote by $\pi : \Au \recht A$ the natural
surjective $^*$-homomorphism, it is clear that $(\pi \ot \io)\al$ will
be a continuous coaction of $(A,\de)$ whenever $\al$ is a continuous
coaction of $(\Au,\deu)$. We briefly discuss the converse: when does
continuous coactions of $(A,\de)$ lift to continuous coactions of
$(\Au,\deu)$?

Fischer has shown (\cite{F}, Proposition 3.26) that a \emph{reduced} continuous coaction has a
unique lift to a continuous coaction of $(\Au,\deu)$. His proof works
for general l.c.\ quantum groups. For regular quantum groups, this is
obvious: a reduced coaction is Morita equivalent to a dual coaction on
a reduced crossed product and it is clear that dual coactions admit a lift.

Another class of coactions for which such a unique lift exists are the
so-called maximal coactions introduced in \cite{F,EKQ}.

\begin{definition} \label{def.maximalcoaction}
Let $(A,\de)$ be a \emph{regular} l.c.\ quantum group.
A continuous coaction $\al : B \recht \M(A \ot B)$ is said to be
\emph{maximal} if the natural surjective $^*$-homomorphism
$$\Aopu \ltimesfull \Ahu \ltimesfull B \recht \K(H) \ot B$$
is an isomorphism.
\end{definition}
Almost by definition, a maximal coaction is a coaction which is Morita
equivalent to a dual coaction on a full crossed product. It is then
clear that maximal coactions admit a unique lift to the universal level.

Up to now we discussed coactions on \cst-algebras. Of course, quantum
groups can coact as well on von Neumann algebras. We refer to \cite{V}
for details.

\begin{definition}
A coaction of a l.c.\ quantum group $(M,\de)$ on a von Neumann algebra
$N$ is a faithful, normal, unital $^*$-homomorphism $\al : N \recht M
\ot N$ satisfying $(\io \ot \al)\al = (\de \ot \io)\al$.

The \emph{crossed product} $\Mh \ltimes N$ is the von Neumann subalgebra of
$\B(H) \ot N$ generated by $\al(N)$ and $\Mh \ot 1$.
\end{definition}

Again the crossed product $\Mh \ltimes N$ carries a \emph{dual
  coaction} of $(\Mh,\dehop)$ and the following biduality result
  holds: $M' \ltimes \Mh \ltimes N \cong \B(H) \ot N$.

\section{Four different pictures of corepresentation theory} \label{sec.fourdiff}

In the preliminary section, we defined unitary corepresentations of
l.c.\ quantum groups and discussed the bijective correspondence with
non-degenerate representations of the universal dual quantum
group. This yields two different pictures of corepresentation
theory. In this section we present two other useful pictures, which
are von Neumann algebraic. These pictures are a major tool in the rest
of the paper.

\subsection{Bicovariant \cst-correspondences}

Before we present the third picture of unitary corepresentation
theory, we give the following definition generalizing the notion of a
correspondence \cite{CJ,Popa,Cnc} from Hilbert spaces to \cst-modules. We
shall need to represent \emph{von Neumann algebras} on \cst-modules.

\begin{definition} \label{def.stricthomo}
Let $N$ be a von Neumann algebra and $\E$ a \cst-$B$-module. A unital
$^*$-homomorphism $\pi : N \recht \cL(\E)$ is said to be \emph{strict}
(or \emph{normal}) if it
is strong$^*$ continuous on the unit ball of $N$.
\end{definition}

Recall that the strong$^*$ topology on the unit ball of $\cL(\E)$
coincides with the strict topology identifying $\cL(\E) =
\M(\K(\E))$. That motivates the terminology of a strict
$^*$-homomorphism.

The following is the almost standard example of a strict
$^*$-homomorphism. Suppose that $N \subset \B(H)$ is a von Neumann
algebra and consider the \cst-$B$-module $H \ot B$. Let $V \in \cL(H
\ot B) = \M( \K(H) \ot B)$ be a unitary operator. Then,
$$N \recht \cL(H \ot B) : x \mapsto V (x \ot 1) V^*$$
is a strict $^*$-homomorphism.

\begin{definition}
Let $M$ and $N$ be von Neumann algebras. We say that a \cst-$B$-module $\E$
is a $B$-corres\-pon\-dence from $N$ to $M$ if we have
\begin{itemize}
\item a strict $^*$-homomorphism $\pil : M \recht \cL(\E)$,
\item a strict $^*$-antihomomorphism $\pir : N \recht \cL(\E)$,
\end{itemize}
such that $\pil(M)$ and $\pir(N)$ commute.
\end{definition}

\begin{notation}
A $B$-correspondence from $N$ to $M$ will be denoted as
$\corresp{\E}{M}{N}$. We will write the left and right module actions
as
$$x \cdot v = \pil(x) v \quad\text{and}\quad v \cdot y = \pir(y) v
\quad \text{for all}\;\; x \in M, y \in N, v \in \E \; .$$
\end{notation}

From \cite{CJ}, we know how to construct a correspondence from
the group von Neumann algebra $\cL(G)$ to $\cL(G)$ given a unitary
representation of $G$. We can do the same thing for l.c.\ quantum
groups.

\begin{proposition}
Let $(A,\de)$ be a l.c.\ quantum group and $X
\in \cL(A \ot \E)$ a unitary corepresentation on a \cst-$B$-module $\E$.

Then, there is a $B$-correspondence $\corresp{H \ot \E}{\Mh}{\Mh}$
given by
$$x \cdot v = X(x \ot 1)X^* v \quad\text{and}\quad v \cdot y = (\Jh
y^* \Jh \ot 1)v \quad\text{for}\quad x,y \in \Mh,\;\; v \in H \ot \E \; .$$
\end{proposition}

To prove this proposition, we only have to observe that $(\io \ot
\pil)(W) = W_{12} X_{13}$, where $\pil : \Mh \recht \cL(\E)$ denotes the
left module action $\pil(x)v = x \cdot v$.
Hence, the left and right module actions commute.

We now want to characterize which $B$-correspondences from $\Mh$ to
$\Mh$ come from a unitary corepresentation.
Let $X \in \cL(A \ot \E)$ be a unitary corepresentation of $(A,\de)$
on a \cst-$B$-module $\E$ and make the $B$-corres\-pon\-dence $\corresp{H
  \ot \E}{\Mh}{\Mh}$. Then, we also have a strict $^*$-homomorphism
$\pi : M' \recht \cL(H \ot \E) : \pi(x) = x \ot 1$ which is covariant
with respect to both the left and right module action of $\Mh$, in the
following precise sense.

\begin{definition}
Let $\corresp{\F}{\Mh}{\Mh}$ be a $B$-corres\-pon\-dence from $\Mh$ to
$\Mh$ and suppose that $\pi : M' \recht \cL(\F)$ is a strict
$^*$-homomorphism.
We say that $\pi$ is \emph{bicovariant} when
$$(\pil \ot \io)\deh(x) = (\pi \ot \io)(\Vh)(\pil(x) \ot 1)(\pi \ot
\io)(\Vh^*)
\quad\text{and}\quad (\pir \ot \Rh)\deh(x) = (\pi \ot
\io)(\Vh^*)(\pir(x) \ot 1)(\pi \ot \io)(\Vh)$$
for all $x \in \Mh$. Here $\pil$ and $\pir$ denote the left and right
module actions of $\Mh$ on $\F$.

In that case, we call $\F$ a \emph{bicovariant $B$-correspondence} and
we write $\bicorresp{\F}{\Mh}{\Mh}{M'}$.
\end{definition}

\begin{remark} \label{rem.sufficient}
Let $\corresp{\F}{\Mh}{\Mh}$ be a $B$-correspondence from $\Mh$ to
$\Mh$. Let $Y \in \cL(\F \ot \Ah)$ be a unitary corepresentation
of $(\Ah,\deh)$, which means that $(\io \ot \deh)(Y) =
Y_{12}Y_{13}$. Suppose that $Y$ is bicovariant in the sense that
\begin{equation}\label{eq.bicovariance}
(\pil \ot \io)\deh(x) = Y(\pil(x) \ot 1)Y^*
\quad\text{and}\quad (\pir \ot \Rh)\deh(x) = Y^*(\pir(x) \ot 1)Y
\end{equation}
for all $x \in \Mh$.

Then, there exists a unique strict $^*$-homomorphism $\pi : M' \recht
\cL(\F)$ which is bicovariant and satisfies $Y = (\pi \ot \io)(\Vh)$.

To show this, you only need $\pil : \Mh \recht \cL(\F)$ satisfying
$(\pil \ot \io)\deh(x) = Y(\pil(x) \ot 1)Y^*$. Writing $X = (\io \ot
\pil)(\Vh)$, this means that $X_{12} \Vh_{13} = Y_{23} X_{12}
Y^*_{23}$.

On the other hand, $Y$ is a corepresentation, which means that $W_{23}
Y_{12} W^*_{23} = Y_{13} Y_{12}$. Recall now that $\Vh = (J\Jh \ot
1)W^*(\Jh J \ot 1)$. So, if we define $\Yh = (J \Jh \ot 1) \Si Y \Si
(\Jh J \ot 1) \in \cL(H \ot \F)$, we get that $\Vh^*_{13} \Yh_{12}
\Vh_{13} = Y_{23} \Yh_{12}$. Together with the formula in the previous
paragraph, we find that
$$(X \Yh)_{12} \Vh_{13} (X \Yh)_{12}^* = Y_{23} \; .$$
From this, it follows that there exists a strict $^*$-homomorphism $\pi : M' \recht
\cL(\F)$ such that $$(X \Yh) (x \ot 1) (X \Yh)^* = 1 \ot \pi(x)$$ for
all $x \in M'$. Then also $Y = (\pi \ot \io)(\Vh)$.

So, a
bicovariant $B$-correspondence is determined by a $B$-correspondence
$\F$ between $\Mh$ and $\Mh$ together with a corepresentation $Y \in
\cL(\F \ot \Ah)$ satisfying the bicovariance relations
\eqref{eq.bicovariance}.
\end{remark}

The following proposition provides the third equivalent picture of
corepresentation theory as \emph{the theory of bicovariant
  $B$-correspondences}.

\begin{proposition}  \label{prop.charbicov}
If $\bicorresp{\F}{\Mh}{\Mh}{M'}$ is a bicovariant $B$-correspondence,
there exists a canonically determined \cst-$B$-module $\E$ and a
corepresentation $X \in \cL(A \ot \E)$ such that
$$\bicorresp{\F}{\Mh}{\Mh}{M'} \;\; \simeq \;\; \bicorresp{H \ot
  \E}{\Mh}{\Mh}{M'}$$
as bicovariant correspondences. So, we get a bijective relation
  between unitary corepresentations on \cst-$B$-modules
  and bicovariant $B$-correspondences.
\end{proposition}

\begin{proof}
Suppose that $\bicorresp{\F}{\Mh}{\Mh}{M'}$ is a bicovariant
$B$-correspondence. Using the technique of the proof of Lemma 5.5 in
\cite{BSV}, we get a strict $^*$-homomorphism $\te : \B(H) \recht
\cL(\F)$ such that $\pi(x) = \te(x)$ for all $x \in M'$ and
$\pir(y) = \te(\Jh y^* \Jh)$ for all $y \in \Mh$. Since $M'$ and
$\Mh'$ generate $\B(H)$ as a von Neumann algebra, the strict
$^*$-homomorphism $\te$ is canonically defined.

Using $\te : \B(H) \recht \cL(\F)$, we get a canonical \cst-$B$-module
$\E$ and an isomorphism $\F \simeq H \ot \E$ such that $\te(x)$
becomes $x \ot 1$ under this isomorphism. One can, for instance,
define $\E$ to be the space of bounded linear maps $v : H \recht \F$
satisfying $v x = \te(x) v$ for all $x \in \B(H)$.

The isomorphism $\F \simeq H \ot \E$ yields $\pil : \Mh \recht \cL(H
\ot \E)$ such that the range of $\pil$ commutes with $\Mh' \ot 1$ and
such that $(\pil \ot \io)\deh(x) = \Vh_{13} (\pil(x) \ot 1) \Vh_{13}^*$.

Write $Z = W^*_{12} (\io \ot \pil)(W) \in \cL(A \ot H \ot \E)$. Then,
$Z$ commutes with $1 \ot \Mh' \ot 1$. On the other hand,
$$\Vh_{24} Z_{123} \Vh_{24}^* = (\io \ot \deh)(W^*)_{124} (\io \ot
\pil \ot \io)(\io \ot \deh)(W) = W^*_{12} W^*_{14} W_{14} (\io \ot
\pil)(W)_{123} = Z_{123} \; .$$
Hence, $Z$ commutes with $1 \ot M' \ot 1$. This implies that there
exists $X \in \cL(A \ot \E)$ such that $Z = X_{13}$. It follows that
$$(\io \ot \pil)(W) = W_{12} X_{13} \; .$$
From this, we conclude that $X$ is a unitary corepresentation of
$(A,\de)$ in $\E$ and that
$$\bicorresp{\F}{\Mh}{\Mh}{M'} \;\; \simeq \;\; \bicorresp{H \ot
  \E}{\Mh}{\Mh}{M'}$$
as bicovariant correspondences.
\end{proof}

\subsection{Bicovariant \wst-bimodules}

We present a fourth picture of corepresentation theory, which only works to describe corepresentations on Hilbert spaces rather then \cst-modules.

\begin{definition}
Let $M,N$ be von Neumann algebras. A \wst-$M$-$N$-bimodule $\E$ is a
\wst-$N$-module equipped with a normal, unital $^*$-homomorphism $\pil
: M \recht \cL(\E)$.
\end{definition}

\begin{proposition}
Let $X \in M \ot \B(K)$ be a unitary corepresentation of a l.c.\
quantum group $(M,\de)$ on a Hilbert space $K$. Consider the
\wst-$\Mh$-module $\Mh \ot K$. Then, there exists a unique normal,
unital $^*$-homomorphism $\pil : \Mh \recht \cL(\Mh \ot K)$ satisfying
$(\io \ot \pil)(W) = W_{12} X_{13}$.

As such, $\Mh \ot K$ becomes a \wst-$\Mh$-$\Mh$-bimodule.
\end{proposition}
\begin{proof}
We define $\pil : \Mh \recht \B(H \ot K) : \pil(x) = X(x \ot
1)X^*$. Then, $(\io \ot \pil)(W) = W_{12} X_{13}$. Hence, $\pil(\Mh)
\subset \Mh \ot K = \cL(\Mh \ot K)$. So, we are done.
\end{proof}

Exactly as we characterized the \cst-correspondences coming from a
corepresentation, we now characterize the \wst-bimodules coming from a
corepresentation. In order to do so, we make use of coactions on
\wst-modules. We refer to the appendix for details on this topic.

If $X \in M \ot \B(K)$ is a unitary corepresentation of $(M,\de)$ on a
Hilbert space $K$, we construct the coaction
$$\ga : \Mh \ot K \recht (\Mh \ot K) \ot \Mh : \ga(z) = (\deh \ot
\io)(z)_{132}$$
on the \wst-$\Mh$-module $\Mh \ot K$
which is compatible with the right coaction $\deh$ on $\Mh$. Moreover,
we observe that
$$\ga \pil = (\pil \ot \io)\deh \; ,$$
where we still write $\ga$ for the associated coaction on $\cL(\Mh \ot
K)$.

This leads to the following definition.
\begin{definition}
Let $(M,\de)$ be a l.c.\ quantum group and let $\F$ be a
\wst-$\Mh$-$\Mh$-bimodule. We say that $\F$ is \emph{bicovariant} if
we have a coaction $\ga : \F \recht \F \ot \Mh$ compatible with the
right coaction $\deh$ on $\Mh$ and satisfying
$$\ga \pil = (\pil \ot \io)\deh \; .$$
\end{definition}

The following result provides the bimodule version of Proposition
\ref{prop.charbicov}.

\begin{proposition}
Let $\F$ be a bicovariant \wst-$\Mh$-$\Mh$-bimodule. Then, there
exists a canonically defined Hilbert space $K$ and a corepresentation
$X \in M \ot \B(K)$ such that $\F \simeq \Mh \ot K$ as bicovariant
\wst-bimodules.
\end{proposition}
\begin{proof}
Define $\F^\ga = \{v \in \F \mid \ga(v) = v \ot 1 \}$. Observe that
$\la v , w \ra$ is invariant under $\deh$ and hence belongs to $\C$
whenever $v,w \in \F^\ga$. So, $\F^\ga$ is a Hilbert space. We shall
show that $\F^\ga$ is the Hilbert space that we are looking for and
that the map $x \ot \xi \mapsto \xi \cdot x$ extends to an isomorphism
from $\Mh \ot \F^\ga$ onto $\F$.

Denote by $N = \cL(\F \oplus \Mh)$ the link algebra and denote by
$\ga$ the right coaction of $(\Mh,\deh)$ on $N$. Define
$$\te : \Mh \recht N : \te(x) = \begin{pmatrix} \pil(x) & 0 \\ 0 & x
\end{pmatrix} \; .$$
Then, $\ga \te = (\te \ot \io)\deh$ by bicovariance of $\E$. From
Proposition 1.22 in \cite{VV}, it follows that $\ga$ is a dual
coaction. This means that the formula
$$N^\ga \recht M \ot N^\ga : z \mapsto (\io \ot \te)(W^*)(1 \ot z)(\io
\ot \te)(W)$$
defines a left coaction of $(M,\de)$ on the fixed point algebra
$$N^\ga = \{ x \in N \mid \ga(x) = x \ot 1 \}$$
such that $N$ is isomorphic with the crossed product $\Mh \ltimes
N^\ga$ and $\ga$ is the dual coaction.

If we now observe that
$$N^\ga = \begin{pmatrix} \cL(\F)^\ga & \F^\ga \\ (\F^\ga)^* & \C
\end{pmatrix}$$
we have found a corepresentation of $(M,\de)$ on the Hilbert space
$\F^\ga$ such that $\F \simeq \Mh \ot \F^\ga$ as bicovariant \wst-bimodules.
\end{proof}

\section{Induction of corepresentations} \label{sec.induction}

We present a new approach to induction of unitary corepresentations of
l.c.\ quantum groups, which
works as well for the induction of corepresentations on
\cst-modules. We first provide some general machinery and start the
induction procedure after Definition \ref{def.coactionwstimprim}.

Let $(M_1,\de_1)$ be a closed quantum subgroup of $(M,\de)$ through
the morphism $(M,\de) \overset{\pi}{\longrightarrow} (M_1,\de_1)$. So,
we have a normal, faithful $^*$-homomorphism $\pih : \Mh_1 \recht \Mh$
satisfying $\deh \pih = (\pih \ot \pih)\deh_1$.

Associated with $\pi$ we have the coaction $\al : M \recht M \ot M_1$
which is formally given by $\al = (\io \ot \pi)\de$ and which, more
precisely, satisfies
$$(\al \ot \io)(W) = W_{13} (\io \ot \pih)(W_1)_{23} \; .$$
Using the modular conjugations, we define as well $\pih' : \Mh_1'
\recht \Mh'$ by
$$\pih' : \Mh_1' \recht \Mh' : \pih'(x) = \Jh \pi(\Jh_1 x \Jh_1) \Jh
\quad\text{for all}\quad x \in \Mh_1' \; .$$

\begin{definition} \label{def.mqhs}
We define $Q = M^\al := \{ x \in M \mid \al(x) = x \ot 1 \}$. The von
Neumann algebra $Q$ should be considered as the \emph{measured quantum
  homogeneous space}.
\end{definition}

Observe that $\de(Q) \subset M \ot Q$. Hence, the restriction of $\de$
to $Q$ defines a left coaction of $(M,\de)$ on $Q$. By definition we
have
$$\Mh \ltimes Q = \bigl( \; \de(Q) \; \cup \; \Mh \ot 1 \; \bigr)\dpr
\; .$$
Observing that $V^* (\Mh \ltimes Q) V = \bigl( \; Q \; \cup \; \Mh \;
\bigr)\dpr \ot 1$, we get that
$$\Mh \ltimes Q \cong \bigl( \; Q \; \cup \; \Mh \;
\bigr)\dpr = (\pih'(\Mh_1'))' \; .$$
We will often identify $\Mh \ltimes Q$ with its image in $\B(H)$.

\begin{definition} \label{def.wstimprim}
Define
$$\cI = \{ v \in \B(H_1,H) \mid v x = \pih'(x) v \;\; \text{for all}\;\; x
\in \Mh_1' \} \; .$$
Defining $\la v,w \ra = v^* w$, the space $\cI$ becomes a
\wst-$\Mh_1$-module. Since $\Mh \ltimes Q = (\pih'(\Mh_1'))'$, we get
that $\cI$ is a \emph{\wst-$(\Mh \ltimes Q)$-$\Mh_1$-imprimitivity bimodule}.
\end{definition}

\begin{remark}
So, we conclude that it is an almost trivial fact that the von Neumann
algebras $\Mh_1$ and $\Mh \ltimes Q$ are Morita equivalent in a von
Neumann algebraic sense.
An important part of the present paper is to define the \emph{locally
  compact quantum homogeneous space $D \subset Q$} such that $\Ah_1$
  is Morita equivalent with $\Ah \ltimesred D$ (and such that $\Ahu_1$
  is Morita equivalent with $\Ahu \ltimesfull D$).
In order to do so, we will make use all the time of the von Neumann
algebraic Morita equivalence $\cI$.
\end{remark}

\begin{definition} \label{def.coactionwstimprim}
Define
$$\al_\cI : \cI \recht \cI \ot \Mh : \al_\cI(v) = \Vh(v \ot 1)(\io \ot
\pih)(\Vh_1^*) \; .$$
Then, $\al_\cI$ is a coaction of $(\Mh,\deh)$ on $\cI$ which is
compatible with the coaction $\deh_1$ on $\Mh_1$. Moreover, if we
equip $\Mh \ltimes Q$ with the dual coaction of $(\Mh,\deh)$, the
right module action of $\Mh \ltimes Q$ on $\cI$ is covariant.
\end{definition}

\emph{We now start the induction procedure.} Let a corepresentation
$X$ of $(M_1,\de_1)$ on a \cst-$B$-module $\E$ be given. So, $X \in
\cL(A_1 \ot \E)$ and $(\de_1 \ot \io)(X) = X_{13} X_{23}$.

Consider the \cst-$B$-module $H \ot \E$. We want to define a strict
$^*$-homomorphism $\pil : \Mh_1 \recht \cL(H \ot \E)$ which is
formally given by the formula $\pil =(\pih \ot \te)\dehopo$, where
$\te : \Ahu_1 \recht \cL(\E)$ is the $^*$-homomorphism associated with
the corepresentation $X$.

\begin{lemma}
There exists a unique strict $^*$-homomorphism $\pil : \Mh_1 \recht
\cL(H \ot \E)$ satisfying
$$(\io \ot \pil)(W_1) = (\io \ot \pih)(W_1)_{12} \; X_{13} \; .$$
\end{lemma}
\begin{proof}
We would like to define $\pil(a) = (\pih \ot \io)(X(a \ot
1)X^*)$. This is somehow delicate, since we would have to give a
meaning to $X(a \ot 1)X^*$ belonging to $\Mh_1 \ot \cL(\E)$ and to
define $(\pih \ot \io)$. We circumvent by defining
$$\pil : \Mh_1 \recht \cL(H \ot \E) : \pil(a) (v \ot 1) \xi = (v \ot
1) X (a \ot 1)X^* \xi$$
for every $\xi \in H_1 \ot \E$ and every $v \in \B(H_1,H)$ satisfying
$v x = \pih(x) v$ for all $x \in \Mh_1$. It is not hard to check that
$\pil(a)$ is a well defined operator in $\cL(H \ot \E)$ and that
$\pil$ is a strict $^*$-homomorphism.
\end{proof}

Equipped with $\pil : \Mh_1 \recht \cL(H \ot \E)$ together with $\pir :
\Mh \recht \cL(H \ot \E) : x \mapsto \Jh x^* \Jh \ot 1$ and $M' \recht
\cL(H \ot \E) : y \mapsto y \ot 1$, we have translated the unitary
corepresentation $X \in \cL(A_1 \ot \E)$ into a \emph{bicovariant
  $B$-correspondence}
\begin{equation}\label{eq.ourcorresp}
\bicorresp{H \ot \E}{\Mh_1}{\Mh}{M'} \; .
\end{equation}
The bicovariance of the above $B$-correspondence can also be expressed
by the coaction
$$\al_{H \ot \E} : H \ot \E \recht \M(H \ot \E \ot \Ah) : \xi \mapsto
\Vh_{13}(\xi \ot 1)$$
which is compatible with the trivial coaction on $B$.

We can now use Definition \ref{def.tensorwst} and Proposition
\ref{prop.productofcoactions} to define
\begin{itemize}
\item the \cst-$B$-module $\Ftil = \cI \rot{\pil} (H \ot \E) \;$;
\item a left module action and a right module action such that we get
  a $B$-correspondence
$\corresp{\Ftil}{\Mh \ltimes Q}{\Mh} \; $;
\item the product coaction $\al_\Ftil$ of $\al_{\cI}$ and $\al_{H \ot \E}$.
\end{itemize}

The product coaction $\al_\Ftil$ is compatible with the trivial coaction
on $B$ and hence yields a corepresentation $Y \in \cL(\Ftil \ot \Ah)$. By
construction (see Remark \ref{rem.sufficient}), we then get the bicovariant $B$-correspondence
$$\bicorresp{\Ftil}{\Mh}{\Mh}{M'} \; .$$
By Proposition \ref{prop.charbicov}, we get a canonically determined
\cst-$B$-module $\Ind \E$ together with a unitary corepresentation
$\Ind X \in \cL(A \ot \E)$ such that
$$\bicorresp{\Ftil}{\Mh}{\Mh}{M'} \;\cong \; \bicorresp{H \ot \Ind
  \E}{\Mh}{\Mh}{M'}$$
as bicovariant correspondences.

Since $Q$ coincides with the fixed point algebra of $\Mh \ltimes Q$
under the dual coaction, we also get a strict $^*$-homomorphism $\rho
: Q \recht \cL(\Ind \E)$ such that
$$\bicorresp{\Ftil}{\Mh \ltimes Q}{\Mh}{M'} \;\cong \; \bicorresp{H \ot \Ind
  \E}{\Mh \ltimes Q}{\Mh}{M'}$$
where, on the right hand side, the left module action $\pil : \Mh \ltimes
  Q \recht \cL(H \ot \Ind \E)$ is determined by
$$(\io \ot \pil)(W) = W_{12} (\Ind X)_{13} \quad\text{and}\quad
\pil(x) = 1 \ot \rho(x) \;\;\text{for all}\; x \in Q \; .$$
By construction we get the following covariance relation:
$$(\io \ot \rho)\de(x) = (\Ind X)^* (1 \ot \rho(x)) (\Ind X)
\quad\text{for all}\;\; x \in \Q \; .$$
Hence, we obtain the expected result that the induced corepresentation
comes with a covariant representation of the measured quantum
homogeneous space.

\begin{definition}
The \cst-$B$-module $\Ind \E$ is called the \emph{induced \cst-$B$-module} of
$\E$ and the unitary corepresentation $\Ind X$ is called the
\emph{induced corepresentation} of $X$.
\end{definition}

Let $B$ and $B_1$ be \cst-algebras. Let $\E$ be a \cst-$B$-module and let $\G$ be a \cst-$B_1$-module. Suppose that $\mu : B \recht \cL(\G)$ is a
non-degenerate $^*$-homomorphism. Then, we have the interior tensor product $\E \rot{\mu} \G$ as a \cst-$B_1$-module. Suppose now that $X \in \cL(A_1
\ot \E)$ is a corepresentation of $(A_1,\de_1)$ on $\E$. Then we have $X \rot{\mu} 1$ as a corepresentation on $\E \rot{\mu} \G$. Since our
construction of the induced corepresentation is completely natural, the following result is obvious.

\begin{proposition} \label{prop.naturality}
We have $\Ind(\E \rot{\mu} \G) \cong \Ind(\E) \rot{\mu} \G$ and
$\Ind(X \rot{\mu} 1) = \Ind(X) \rot{\mu} 1$ in a natural way.

The
representation of $Q$ on $\Ind(\E \rot{\mu} \G)$ is intertwined with
the representation $Q \rot{\mu} 1$ on $\Ind(\E) \rot{\mu} \G$.
\end{proposition}

\section{First imprimitivity theorem}  \label{sec.firstimprim}

In the previous section we defined the induced corepresentation $\Ind
X$ of a corepresentation $X$ of a closed quantum subgroup
$(M_1,\de_1)$ of $(M,\de)$. Such an induced corepresentation comes
with a covariant representation of the measured quantum homogeneous
space $Q \recht \cL(\Ind \E)$.

A natural question is now of course if an imprimitivity result holds.
More precisely, let $Z \in \cL(A \ot \F)$ be a corepresentation of
$(A,\de)$ on a \cst-$B$-module $\F$ and $\rho : Q \recht \cL(\F)$ a
strict $^*$-homomorphism which is covariant. Does there exist a
corepresentation $X$ of $(A_1,\de_1)$ on a \cst-$B$-module $\E$ such
that
$$(\F \; , \; Z \; ,\; \text{rep.\ of}\; Q) \; \cong \;  (\Ind \E \; , \; \Ind
X \; , \;  \text{rep.\ of}\;
Q) \; ?$$

It is quite clear that the answer is negative in general. If $M_1 =
\C$, the one-point subgroup, it is obvious to check that $Q = M$ and
that the induced corepresentations are the multiples of the regular
corepresentation $W$ of $(M,\de)$. The question in the previous
paragraph becomes the following: is every pair $(X,\rho)$ of a
corepresentation $X$ of $(M,\de)$ and a covariant representation
$\rho$ of $M$ isomorphic with a multiple of the regular
corepresentation of $(M,\de)$ and the standard representation of $M$
on $H$? But, this question is equivalent with the question
$$\Ahu \ltimesfull A \cong \K(H) \;\; ?$$
This property is precisely the \emph{strong regularity} of the quantum
group $(A,\de)$, see Definition \ref{def.strongreg}.

Conclusion: we can only hope for an imprimitivity result if strong
regularity holds. Otherwise, imprimitivity already fails for the
one-point subgroup.

\begin{theorem}[First Imprimitivity Theorem] \label{thm.firstimprim}
Let $(M,\de)$ be a strongly regular locally compact quantum group. Let
$(M_1,\de_1)$ be a closed quantum subgroup and $Q \subset M$ the
measured quantum homogeneous space in the sense of Definition
\ref{def.mqhs}.

A corepresentation $Z \in \cL(A \ot \cF)$ of $(M,\de)$ on the
\cst-$B$-module $\cF$ is induced from
a corepresentation of $(M_1,\de_1)$ if and only if there exists a
strict $^*$-homomorphism $\rho : Q \recht \cL(\cF)$ such that
\begin{equation}\label{eq.covariance}
(\io \ot \rho)\de(x) = Z^* (1 \ot \rho(x)) Z \quad\text{for all}\;\;
x \in Q \; .
\end{equation}
\end{theorem}
\begin{proof}
It is clear that we only have to prove one implication. So, let
$Z \in \cL(A \ot \cF)$ be a corepresentation of $(M,\de)$ on the
\cst-$B$-module $\cF$ and let $\rho : Q \recht \cL(\cF)$ be a strict
$^*$-homomorphism satisfying the covariance relation
\eqref{eq.covariance}.

To obtain the \cst-$B$-module $\E$ and a corepresentation of
$(M_1,\de_1)$ on it, we perform exactly the inverse of the induction
procedure, tensoring with the inverse of the \wst-imprimitivity
bimodule $\cI$ defined in Definition \ref{def.wstimprim}.

We claim that there exists a unique normal $^*$-homomorphism $\pil : \Mh
\ltimes Q \recht \cL(H \ot \F)$ such that
$$(\io \ot \pil)(W) = W_{12} Z_{13} \quad\text{and}\quad \pil(x) = 1
\ot \rho(x) \;\;\text{for all}\;\; x \in Q \; .$$
Indeed, it suffices to define $\pil(z) = Z (\io \ot \rho)(V(z \ot
1)V^*) Z^*$. For all $z \in \Mh \ltimes Q$, the element $V(z
\ot 1)V^*$ belongs to $\B(H) \ot Q$. From Lemma \ref{lem.extendstrict}
we know that we can extend $\io \ot \rho$ to the von Neumann algebra
$\B(H) \ot Q$. It is easy to check that $\pil$, once well-defined,
satisfies the required conditions.

Using the anti-homomorphism $\pir : \Mh \recht \cL(H \ot \F) : \pir(x)
= \Jh x^* \Jh \ot 1$ and the homomorphism $M' \recht \cL(H \ot \F) : y
\mapsto y \ot 1$, we get a bicovariant $B$- correspondence
$$\bicorresp{H \ot \F}{\Mh \ltimes Q}{\Mh}{M'} \; .$$

In Definition \ref{def.wstimprim} we defined the \wst-imprimitivity
bimodule $\cI$. We can define its inverse as
$$\cI^* = \{ v \in \B(H,H_1) \mid x v = v \pih'(x) \quad\text{for
  all}\;\; x \in \Mh_1' \} \; .$$
Then, $\cI^*$ is a \wst-$\Mh_1$-$\Mh \ltimes Q$-bimodule. We can again
  define, as in Definition \ref{def.coactionwstimprim} a coaction of
  $(\Mh,\deh)$ on $\cI^*$.

We can then define, using Definition \ref{def.tensorwst} and
Proposition \ref{prop.productofcoactions},
\begin{itemize}
\item the \cst-$B$-module $\Etil = \cI^* \rot{\pil} (H \ot \cF)$,
\item a left and a right module action such that we get a
  $B$-correspondence
$\corresp{\Etil}{\Mh_1}{\Mh}\; ,$
\item the product coaction $\al_\Etil$ of $(\Mh,\deh)$ on $\Etil$.
\end{itemize}
Hence, we have a bicovariant $B$-correspondence
$$\bicorresp{\Etil}{\Mh_1}{\Mh}{M'} \; .$$
The homomorphism $M' \recht \cL(\Etil)$ and the anti-homomorphism $\Mh
\recht \cL(\Etil)$ are covariant. From the strong regularity of
$(M,\de)$ it follows that we find a canonically determined
\cst-$B$-module $\E$ such that $\Etil \cong H \ot \E$ where the
isomorphism intertwines the homomorphism $M' \recht \cL(\Etil)$ with
$x \mapsto x \ot 1$ and the anti-homomorphism $\Mh
\recht \cL(\Etil)$ with $y \mapsto \Jh y^* \Jh \ot 1$.

Exactly as in the proof of Proposition \ref{prop.charbicov} the
homomorphism $\Mh_1 \recht \cL(\Etil)$ is intertwined with a
homomorphism $\pil : \Mh_1 \recht \cL(H \ot \E)$ such that
$$(\io \ot \pil)(W_1) = (\io \ot \pih)(W_1)_{12} X_{13}$$
where $X \in \cL(A_1 \ot \E)$ is a corepresentation of $(M_1,\de_1)$
on the $B$-module $\E$. We get
$$\bicorresp{\Etil}{\Mh_1}{\Mh}{M'} \; \;\cong \; \; \bicorresp{H \ot
  \E}{\Mh_1}{\Mh}{M'}$$
as bicovariant correspondences.

It remains to prove that $\F = \Ind \E$ and $Z = \Ind X$. For this it
suffices to observe that the interior tensor product of $\cI$ and
$\cI^*$ is canonically isomorphic with $\Mh \ltimes Q$ as a \wst-$\Mh
\ltimes Q$-bimodule equipped with the dual coaction. Hence,
$$H \ot \Ind \E = \cI \rot{\Mh_1} (H \ot \E) = \cI \rot{\Mh_1} \cI^* \rot{\Mh
  \ltimes Q} (H \ot \cF) = (\Mh \ltimes Q) \rot{\Mh \ltimes Q} (H \ot
  \cF) = H \ot \cF \; .$$
\end{proof}

\section{Quantum homogeneous spaces and Mackey imprimitivity} \label{sec.qhs}

We fix a locally compact quantum group $(M,\de)$. {\it We suppose
  throughout this section that $(M,\de)$ is strongly regular.}
We fix a closed quantum
subgroup $(M_1,\de_1)$. Recall that $Q \subset M$ denotes the measured
quantum homogeneous space.

We shall prove the following crucial results.

\begin{theorem} \label{thm.existqhs}
There exists a unique \cst-subalgebra $D \subset Q$ satisfying
\begin{itemize}
\item $D$ is strongly dense in $Q$,
\item $\de(D) \subset \M(A \ot D)$ and $\de : D \recht \M(A \ot D)$ is
  a continuous coaction of $(A,\de)$ on $D$,
\item $\de(Q) \subset \M(\K(H) \ot D)$ and the $^*$-homomorphism $\de
  : Q \recht \cL(H \ot D)$ is strict.
\end{itemize}
We call $D$ the \emph{quantum homogeneous space}.
\end{theorem}

\begin{theorem} \label{thm.morita}
There exist canonical covariant Morita equivalences
$$\Ah \ltimesred D \Morita \Ah_1
\quad\text{and}\quad \Ahu \ltimesfull D \Morita \Ahu_1 \; .$$
\end{theorem}

As we shall see, the uniqueness statement in Theorem
\ref{thm.existqhs} is not
so hard to prove and valid without the assumption on strong
regularity. Our existence proof of $D$ uses the strong regularity
assumption (in fact, regularity suffices) but it is not excluded that
$D$ even exists without regularity assumptions. We recall however that
we cannot hope for an imprimitivity theorem in the non-regular case.

\begin{remark}
The precise meaning of the statement Theorem \ref{thm.morita} is the following. There exist a \cst-$\Ahu_1$-module $\Efull$, a natural isomorphism
$\K(\Efull) \cong \Ahu \ltimesfull D$ and a right coaction of $(\Ah,\deh)$ on $\Efull$ that is compatible with the coaction $(\io \ot \pih)\dehu_1$
on $\Ahu_1$ and coincides with the dual coaction on $\Ahu \ltimesfull D$. Composing the Morita equivalence $\Efull$ between $\Ahu \ltimesfull D$ and
$\Ahu_1$ on one side with the regular representation, we get the Morita equivalence $\Ah \ltimesred D \Morita \Ah_1$.
\end{remark}

The right coaction $\al : M \recht M \ot M_1$ of $(M_1,\de_1)$ on $M$ by right translation (see the beginning of Section \ref{sec.induction}),
restricts to a continuous right coaction of $(A_1,\de_1)$ on $A$. The reduced crossed product $A \rtimesred \Ahop_1$ is given as the closed linear
span of $\al(A) (1 \ot \Jh_1 \Ah_1 \Jh_1)$.

\begin{corollary} \label{cor.morita}
There is a natural covariant Morita equivalence
$$D \Morita A \rtimesred \Ahop_1 \; .$$
\end{corollary}

\begin{proof}
Observe that we have covariant isomorphisms
$$A \rtimesred \Ahop_1 \cong [A \; \Jh \pih(\Ah_1) \Jh]
\cong [JAJ \; \pih(\Ah_1)] \cong \Ah_1 \rtimesred \Aop \; ,$$ where we consider the crossed product of $\Ah_1$ equipped with the right coaction $(\io
\ot \pih)\deh_1$ of $(\Ah,\deh)$.

Since we have a covariant Morita equivalence $\Ah \ltimesred D
\Morita \Ah_1$, we can take the crossed product
in order to obtain a covariant Morita equivalence
$$(\Ah \ltimesred D) \rtimesred \Aop
\Morita \Ah_1 \rtimesred \Aop \; .$$ The biduality theorem gives a natural covariant Morita equivalence
$$(\Ah \ltimesred D) \rtimesred \Aop
\Morita D$$
and then we are done.
\end{proof}

\begin{remark}
We know of course that whenever $G_1$ is a closed subgroup of a
locally compact group $G$, the action of $G_1$ on $G$ by right
translations is proper. Hence, the full and reduced crossed products
$C_0(G) \rtimesfull G_1$ and $C_0(G) \rtimesred G_1$ coincide.

There are strong indications that the same result is no longer valid in general, even for strongly regular l.c.\ quantum groups. Nevertheless, if
either $(A,\de)$ is co-amenable (which means that $\Au = A$) or $(A_1,\de_1)$ is amenable, the full and reduced crossed products $A \rtimesfull
\Ahopu_1$ and $A \rtimesred \Ahop_1$ coincide. The second part is of course obvious. So, suppose that $(A,\de)$ is co-amenable. Observe that this is
for instance the case when $A = C_0(G)$ and hence, this case covers the group case of the previous paragraph. Using twice the co-amenability of
$(A,\de)$, we have
$$A \rtimesfull \Ahopu_1 \cong \Au \rtimesfull \Ahopu_1
\cong \Ahu_1 \rtimesfull \Aopu \cong \Ah_1 \rtimesfull \Aopu \cong \Ah_1 \rtimesred \Aop \cong A \rtimesred \Ahop_1 \; .$$ In general, although I do
not have an explicit example, it might very well be that $(A,\al) \cong (A_1 \ltimesred B, \hat{\be})$, where $\be : B \recht \M(\Ah_1 \ot B)$ is a
sufficiently non-trivial reduced coaction and $\hat{\be}$ is the dual right coaction of $(A_1,\de_1)$ on the crossed product. The statement $A
\rtimesfull \Ahopu_1 = A \rtimesred \Ahop_1$ comes down to saying that $\be$ is as well a maximal coaction. There seems to be no reason why this
should always be the case in a non-amenable, non-co-amenable situation.
\end{remark}

In order to construct $D$ we have to look at a covariant induction
procedure. This has an independent interest. Indeed, we shall not only
prove the Morita equivalence
$$\Ahu \ltimesfull D \Morita \Ahu_1$$
but we also want that this Morita equivalence is covariant for a natural coaction of $(\Ah,\deh)$, compatible with the dual coaction on $\Ahu
\ltimesfull D$ and the comultiplication on $\Ahu_1$.

Let $\E$ be a \cst-$B$-module. Suppose that we have a
coaction $\be : \E \recht \M(\E \ot \Ah)$ compatible with a continuous
coaction on $B$ that we also denote by $\be : B \recht \M(B \ot \Ah)$.

Let $X \in \cL(A_1 \ot \E)$ be a corepresentation of $(M_1,\de_1)$ on
$\E$ satisfying the following compatibility relation with $\be$:
\begin{equation}\label{eq.nogcovar}
(\io \ot \be)\bigl( X (1 \ot v) \bigr) = (\io \ot \pih)(W_1)_{13}
X_{12} \bigl( 1 \ot \be(v)  \bigr) \quad\text{for all}\;\; v \in \E \; .
\end{equation}
If we consider the $^*$-homomorphism $\te : \Ahu_1 \recht \cL(\E)$
associated with the corepresentation $X$, Equation \eqref{eq.nogcovar}
becomes $\be(\te(a) v) = (\te \ot \pih)\deh_1(a) \be(v)$ for all $a
\in \Ahu_1$ and $v \in \E$.

Let $\F = \Ind \E$ be the induced \cst-$B$-module with induced
corepresentation $Y=\Ind X$ of $(A,\de)$ on $\F$.
We shall construct an induced coaction $\Ind \be$ on the induced
\cst-module $\Ind \E$.

Recall that, by
definition of $\Ind \E$, we have
an isomorphism
$$\Phi : \cI \rot{\pil} (H \ot \E) \recht H \ot \F$$
where $\pil : \Mh_1 \recht \cL(H \ot \E)$ is determined by $(\io \ot
\pil)(W_1) = (\io \ot \pih)(W_1)_{12} X_{13}$.

\begin{proposition} \label{prop.inducedcoaction}
There exists a unique coaction $\Ind \be : \Ind \E \recht
\M(\Ind \E \ot \Ah)$ of $(\Ah,\deh)$ on the
induced \cst-$B$-module $\Ind \E$ which is compatible with the
coaction $\be$ on $B$ and satisfies
\begin{equation}\label{eq.ourcoaction}
(\io \ot \Ind \be) \Phi(v \rot{\pil} x) = W_{13} (\Phi \ot
\io)\bigl( (v \ot 1) \rot{\pil \ot \io} W^*_{13} (\io \ot \be)(x)
\bigr) \quad\text{for all}\;\; v \in \cI, x \in H \ot \E \; .
\end{equation}
Writing $\ga$ for $\Ind \be$, writing $\ga$ as well for the
associated coaction on $\K(\Ind \E)$ and writing $Y = \Ind X$, we get
that
\begin{equation}\label{eq.nogeen}
Q \subset \cL(\Ind \E)^\ga \quad\text{and}\quad (\io \ot \ga)(Y) =
W_{13} Y_{12} \; .
\end{equation}
\end{proposition}
Remark that the defining Equation \eqref{eq.ourcoaction} makes sense
as an equality in $\M(H \ot \E \ot H)$.

\begin{proof}
Write $\Etil = H \ot \E$.
Define the coaction $\eta : \Etil \recht \M(\Etil \ot \Ah)$ by
$\eta(v) = W^*_{13} (\io \ot \be)(v)$ and observe that $\eta$ is
compatible with the coaction $\be$ on $B$. We still write $\eta$ for the
associated coaction on $\K(\Etil)$. It is easily verified that
$\eta(\pil(x)) = \pil(x) \ot 1$ for all $x \in \Mh_1$.

If we equip $\cI$ and $\Mh_1$ with the trivial coaction of $(\Mh,\deh)$, the
homomorphism $\pil : \Mh_1 \recht \cL(\Etil)$ is covariant in the
sense of Definition \ref{def.covariant}. Hence, Proposition
\ref{prop.productofcoactions} yields a product coaction $\eta_1$ of
$(\Ah,\deh)$ on
$\Ftil = \cI \rot{\pil} \Etil$ such that
$$\eta_1(v \rot{\pil} x) = (v \ot 1) \rot{\pil \ot \io} W^*_{13} (\io
\ot \be)(x)$$
for all $v \in \cI$ and $x \in \Etil$.

We have an isomorphism $\Phi : \Ftil \recht H \ot \F$ and this allows
to define $\eta_2 : H \ot \F \recht \M(H \ot \F \ot \Ah)$ such that
$$
\eta_2(\Phi(x)) = W_{13} (\Phi \ot \io)(\eta_1(x)) \quad\text{for
  all}\;\; x \in \Ftil \; .$$

We claim that $\eta_2$ is invariant under the right module action of
$\Mh$ on $H \ot \F$. Indeed, for $a \in \Mh$ and $x \in \Ftil$, we have
$$\eta_2 \bigl( (\Jh a^* \Jh \ot 1) \Phi(x) \bigr)
=\eta_2\bigl( \Phi(\pir(a) x) \bigr) = W_{13} (\Phi \ot
\io)\bigl(\eta_1(\pir(a)x)\bigr) \; .$$
We then observe that for $a \in \Mh$, $v \in \cI$ and $y \in \Etil$
$$\eta_1(\pir(a) (v \rot{\pil} y)) = (v \ot 1) \rot{\pil \ot \io}
W^*_{13} (\io \ot \be)( (\Jh a^* \Jh \ot 1) y)
= (v \ot 1) \rot{\pil \ot \io} ((\Jh \ot J)\dehop(a^*)(\Jh \ot
J))_{13} W^*_{13} (\io \ot \be)(y) \; .$$
We conclude that
$$W_{13} (\Phi \ot
\io)\bigl(\eta_1(\pir(a)x)\bigr) = W_{13} ((\Jh \ot J)\dehop(a^*)(\Jh \ot
J))_{13} (\Phi \ot \io)\eta_1(x) = (\Jh a^* \Jh \ot 1 \ot 1)
\eta_2(\Phi(x)) \; .$$
This proves our claim.

On the other hand, in the
  construction of the induced module $\Ind \E$, we used the product
  coaction $\al_\Ftil$ of $\al_\cI$ and $\al_{H \ot \E}$. It is clear
  that the coactions $\al_\Ftil$ and $\eta_1$ commute. This implies
  that $\eta_2$ is invariant under the representation $M' \ot 1$ of
  $M'$ on $H \ot \F$.

So, we have shown that $\eta_2$ is invariant under $M' \ot 1$ as well
as $\Mh' \ot 1$. Hence, there exists a non-degenerate linear map $\ga
: \F \recht \M(\F \ot \Ah)$ such that $\eta_2 = \io \ot \ga$. Since
$\eta_1$ is a coaction, the map $x \mapsto W^*_{13} (\io \ot \ga)(x)$
defines a coaction of $(\Ah,\deh)$ on $H \ot \F$. This implies that
$\ga$ is as well a coaction of $(\Ah,\deh)$ on $\F$. We define $\Ind \be := \ga$.

By definition $\eta_1$ is invariant under the left module action of
$\Mh \ltimes Q$ on $\Ftil$. It is then clear that $\ga$ satisfies
\eqref{eq.nogeen}.
\end{proof}

We are now ready to prove the main Theorems \ref{thm.existqhs} and \ref{thm.morita}

\begin{proof}[Proof of Theorem \ref{thm.existqhs}]
Consider the \cst-$\Ah_1$-module $\Ah_1$ equipped with the regular
corepresentation $W_1 \in \M(A_1 \ot \Ah_1)$ and the coaction $\be := (\io
\ot \pih)\deh_1$ of $(\Ah,\deh)$.

Define the \cst-$\Ah_1$-module $\J = \Ind \Ah_1$ together with the
induced corepresentation $X \in \M(A \ot \J)$, the strict
$^*$-homomorphism $\te : Q \recht \cL(\J)$ and the induced coaction
$\ga = \Ind \be$ of $(\Ah,\deh)$ on $\J$ as in Proposition
\ref{prop.inducedcoaction}.

Continue writing $\ga$ for the coaction of $(\Ah,\deh)$ on
$\K(\J)$. Then, we have a strict $^*$-homomorphism $\te : Q \recht
\cL(\J)$ such that $\ga(\te(x)) = \te(x) \ot 1$ for all $x \in
Q$. We also have $(\io \ot \ga)(X) = W_{13} X_{12}$.

We claim that $\te : Q \recht \cL(\J)^\ga$ is a $^*$-isomorphism. Using
the regular representation $\Ah_1 \recht \B(H_1)$, we get, using
Proposition \ref{prop.naturality} that
$$\J \rot{\Ah_1} H_1 = \Ind(\Ah_1) \rot{\Ah_1} H_1 = \Ind(H_1) = H \;
.$$
It is straightforward to check that, under these identifications,
$$\te(x) \rot{\Ah_1} 1 = x \quad\text{for all}\;\; x \in Q$$
and $X \rot{\Ah_1} 1 = W$. So, we get an injective $^*$-homomorphism $\K(\J)
\recht \Mh \ltimes Q$ which intertwines the coaction $\ga$ with the
dual coaction on $\Mh \ltimes Q$. Since the fixed point algebra of
$\Mh \ltimes Q$ under the dual coaction is precisely $Q$, we have
proved our claim.

Combining Theorem \ref{thm.quantumlandstad} and Remark
\ref{rem.strictcontinuity} with the isomorphism $\cL(\J)^\ga \cong Q$, we conclude
that there exists a strongly dense \cst-subalgebra $D \subset Q = \cL(\J)^\ga$ such that
\begin{itemize}
\item $\de : D \recht \M(A \ot D)$ is a continuous coaction of
  $(A,\de)$ on $D$;
\item $\de : Q \recht \M(\K(H) \ot D)$ is well defined and strict.
\end{itemize}
So, we have proved the existence part of Theorem \ref{thm.existqhs}.

To prove uniqueness, suppose that $D_1$ and $D_2$ satisfy the
conditions in the theorem. Then, we get that $D_1 = [(\om \ot \io)\de(D_1) \mid
\om \in \B(H)_*]$ by continuity of the coaction. Using the strictness
of $\de : Q \recht \M(\K(H) \ot D_1)$ and the fact that $D_1$ as well
as $D_2$ are dense
in $Q$, we obtain that
\begin{align*}
D_1 &= [D_1 D_1] =  [(\om \ot \io)\bigl(\de(D_1) (\K(H) \ot D_1)\bigr) \mid
\om \in \B(H)_*] = [(\om \ot \io)\bigl(\de(Q) (\K(H) \ot D_1)\bigr) \mid
\om \in \B(H)_*] \\ &= [(\om \ot \io)\bigl(\de(D_2) (\K(H) \ot D_1)\bigr) \mid
\om \in \B(H)_*] = [D_2 D_1] \; .
\end{align*}
By symmetry, we find that $D_2 = [D_1 D_2]$. Taking the adjoint, this
gives $D_2 = [D_2 D_1]$ and we conclude that $D_1 = D_2$.
\end{proof}

\begin{proof}[Proof of Theorem \ref{thm.morita}]
The statement $\Ah \ltimesred D \Morita \Ah_1$
has already been shown in the proof of Theorem \ref{thm.existqhs}.

Loosely speaking, the first imprimitivity theorem
\ref{thm.firstimprim} tells us that representations of $\Ahu_1$ on
\cst-$B$-modules are in one to one correspondence to covariant
representations of the pair $(\Ahu,Q)$ on \cst-$B$-modules. In order
to prove the theorem, it suffices to show that
\begin{equation}\label{eq.statement}
\text{any rep.\ of $\Ahu \ltimesfull D$ on a \cst-$B$-module $\F$ extends to a strict $^*$-homomorphism $Q \recht \cL(\F)$.}
\end{equation}
We shall show this
statement at the end of the proof.

Consider the \cst-$\Ahu_1$-module $\Ahu_1$ equipped with the universal
corepresentation $\Wu_1 \in \M(A_1 \ot \Ahu_1)$ and the coaction $\be
:= (\io
\ot \pih)\dehu_1$ of $(\Ah,\deh)$.

Define the \cst-$\Ahu_1$-module $\Efull = \Ind \Ahu_1$ together with the induced corepresentation $X \in \M(A \ot \Efull)$, the strict
$^*$-homomorphism $\te : Q \recht \cL(\Efull)$ and the induced coaction $\ga = \Ind \be$ of $(\Ah,\deh)$ on $\Efull$ as in Proposition
\ref{prop.inducedcoaction}. By the covariance of $X$ and $\te$, we have an associated representation $\Ahu \ltimesfull D \recht \cL(\Efull)$, which
intertwines the dual coaction on $\Ahu \ltimesfull D$ with the coaction $\ga$ on $\Efull$. We claim that this isomorphism identifies $\Ahu
\ltimesfull D$ with $\K(\Efull)$.

Consider the \cst-$(\Ahu \ltimesfull D)$-module $\Ahu \ltimesfull D$. Using \eqref{eq.statement}, we get a covariant representation of $\Ahu$ and $Q$
in $\M(\Ahu \ltimesfull D) = \cL(\Ahu \ltimesfull D)$. By the first imprimitivity theorem \ref{thm.firstimprim}, we get a \cst-$(\Ahu \ltimesfull
D)$-module $\Eprimefull$ together with a representation of $\Ahu_1$ on $\Eprimefull$ such that $\Ahu \ltimesfull D = \Ind \Eprimefull$. Then, by
Proposition \ref{prop.naturality}, we get
$$\Efull \rot{\Ahu_1} \Eprimefull = (\Ind \Ahu_1) \rot{\Ahu_1}
\Eprimefull = \Ind (\Ahu_1 \rot{\Ahu_1} \Eprimefull) = \Ind \Eprimefull = \Ahu \ltimesfull D \; .$$ Conversely,
$$\Ind(\Eprimefull \rot{\Ah \ltimesfullsmall D} \Efull) = (\Ind
\Eprimefull) \rot{\Ah \ltimesfullsmall D} \Efull = (\Ahu \ltimesfull D) \rot{\Ah \ltimesfullsmall D} \Efull = \Efull = \Ind(\Ahu_1) \; .$$ From this
it follows that $\Eprimefull \rot{\Ah \ltimesfullsmall D} \Efull = \Ahu_1$. Since we have found the inverse module $\Eprimefull$, our claim is
proved.

It remains to prove \eqref{eq.statement}. So, let $Z \in
\cL(A \ot \F)$ be a corepresentation of $(A,\de)$ on a \cst-$B$-module
$\F$ and let $\rho : D \recht \cL(F)$ be a non-degenerate
$^*$-homomorphism which is covariant in the sense that
$$(\io \ot \rho)\de(x) = Z^* (1 \ot \rho(x)) Z \quad\text{for all}\;\;
x \in D \; .$$
Since $\de : Q \recht \M(\K(H) \ot D)$ is strict, we can define a
strict $^*$-homomorphism
$$\mu : Q \recht \cL(H \ot \F) : \mu(x) = Z (\io \ot \rho)\de(x) Z^*
\; .$$
But then, $\mu(x) = 1 \ot \rho(x)$ for all $x \in D$. Since $D$ is
dense in $Q$ and since $\mu$ is strict, it follows that there exists a
strict $^*$-homomorphism $\eta : Q \recht \cL(\F)$ such that $\mu(z) =
1 \ot \eta(z)$ for all $z \in Q$. Then, $\eta$ is the extension of
$\rho$ that we were looking for.
\end{proof}

The major tool used in the proof of Theorem \ref{thm.existqhs} is the characterization of reduced crossed products. This is a quantum version of a
theorem of Landstad (Theorem 3 in \cite{Lan}).

\begin{theorem} \label{thm.quantumlandstad}
Let $(A,\de)$ be a locally compact quantum group. Let $\be : B \recht
\M(B \ot \Ah)$ be a reduced continuous coaction of $(\Ah,\deh)$ on the
\cst-algebra $B$. Then, the following conditions are equivalent.
\begin{enumerate}
\item There exists a \cst-algebra $D$ and a continuous coaction $\al :
  D \recht \M(A \ot D)$ such that $(B,\be) \cong (\Ah \ltimesred
  D,\alh)$, where $\alh$ denotes the dual coaction.
\item There exists a corepresentation $X \in \M(A \ot B)$ of $(A,\de)$
  in $B$ which is covariant in the sense that
$$(\io \ot \be)(X) = W_{13} X_{12} \; .$$
\end{enumerate}
If the second condition is fulfilled, $D$ can be taken as the unique
\cst-subalgebra of $\M(B)^\be$ satisfying
\begin{itemize}
\item the map $\al : x \mapsto X^*(1 \ot x)X$ defines a continuous coaction
  of $(A,\de)$ on $D$;
\item $B$ is the closed linear span of $\{ x (\om \ot \io)(X) \mid x
  \in D, \om \in \B(H)_* \}$.
\end{itemize}
Moreover, an explicit $^*$-isomorphism $B \recht \Ah \ltimesred D$ is then
given by $\eta : z \mapsto X^* \be(z)_{21} X$.
\end{theorem}

\begin{proof}
In the course of this proof, we denote by $[X]$ the closed linear span of a subset of a \cst-algebra. Suppose first that $\al : D \recht \M(A \ot D)$
is a continuous coaction of $(A,\de)$ on the \cst-algebra $D$. Then, the crossed product $\Ah \ltimesred D$ is defined as $[\al(D)(\Ah \ot 1)]$. It
is clear that we can take $X = W \ot 1 \in \M(A \ot (\Ah \ltimesred D))$.

Suppose next that the second condition is fulfilled. If we write $\eta
: B \recht \M(\K \ot B) : \eta(z) = X^*\be(z)_{21} X$, we observe that
$(\io \ot \eta)(X) = W \ot 1$. Hence, there is a uniquely defined
faithful non-degenerate $^*$-homomorphism $\te : \Ah \recht \M(B)$
such that $X = (\io \ot \te)(W)$. Moreover $\be \te = (\te \ot
\io)\deh$.

We first prove the uniqueness statement. Suppose that $D \subset
\M(B)^\be$ is a \cst-algebra that satisfies both conditions in the
theorem. Then,
\begin{align*}
[(\om \ot \io)\eta(z) \mid z \in B, \om \in \B(H)_*] &=
[(\om \ot \io)\eta(x \te(a)) \mid x \in D, a \in \Ah, \om \in \B(H)_*]
\\ &=[(\om \ot \io)\bigl( X^*(1 \ot x) X (a \ot 1) \bigr) \mid x \in D, a
\in \Ah, \om \in \B(H)_*]
\\ &=[(\om \ot \io)\bigl( X^*(1 \ot x) X \bigr) \mid x \in D, \om \in \B(H)_*] = D \; .
\end{align*}
Since the left hand side does not depend on $D$, uniqueness of $D$ has
been proved.

In order to prove existence of $D$, define
$$D := [(\om \ot \io)\eta(z) \mid z \in B, \om \in \B(H)_*] \; .$$
We first show that $D$ is a \cst-algebra. Since the coaction $\be$ of
$(\Ah,\deh)$ on $B$ is continuous, we get that $[\be(B) (1 \ot JAJ)] =
[(1 \ot JAJ) \be(B)]$, since this space is exactly the crossed product
of $B$ with the coaction $\be$. By regularity of $(A,\de)$, we also
know that $\K(H) = [JAJ \; \Ah]$. So,
\begin{align*}
[\eta(B) (\K(H) \ot 1) \eta(B)] &= [\eta(B) (JAJ \; \Ah \ot 1)
\eta(B)] = [X^* \bigl( \be(B) (1 \ot JAJ) \bigr)_{21} \; X (\Ah \ot 1)
X^* \; \be(B)_{21} X ] \\ &= [(JAJ \ot 1) X^* \bigl( \be(B) \;
\be(\te(\Ah)) \; \be(B)\bigr)_{21} X] = [(JAJ \ot 1) \eta(B)] \; .
\end{align*}
Applying $(\om \ot \io)$ on both sides of this equality, we obtain
that $D = [DD]$. Hence, $D$ is a \cst-subalgebra of $\M(B)^\be$.

Define $\al : D \recht \M(A \ot B) : \al(x) = X^*(1 \ot x)X$. Then,
\begin{align*}
[\al(D) (A \ot 1)] &= [(\om \ot \io \ot \io)(X^*_{23} X^*_{13}
\be(B)_{31} X_{13} X_{23})(A \ot 1) \mid \om \in \B(H)_* ]
\\ &=[(\om \ot \io \ot \io)\bigl( V_{12} X_{13}^* V_{12}^* \; \be(B)_{31}
\; V_{12} X_{13} V_{12}^*\bigr)(A \ot 1) \mid \om \in \B(H)_* ]
\\ &=[(\om \ot \io \ot \io)\bigl( V_{12} \eta(B)_{13} V_{12}^* (\K(H) \ot
A \ot 1) \bigr) \mid \om \in \B(H)_* ]
\\ &=[(\om \ot \io \ot \io)\bigl( V_{12} \eta(B)_{13} (1 \ot A \ot 1)
\mid \om \in \B(H)_* ]
\end{align*}
because $V \in \M(\K(H) \ot A)$. From the regularity of $V$, it
follows that $[(\K(H) \ot 1) V (1 \ot A)] = \K(H) \ot A$ and hence,
$$[\al(D) (A \ot 1)] = A \ot [(\om \ot \io)\eta(B) \mid \om \in
\B(H)_* ] = A \ot D \; .$$ So, $\al$ defines a continuous coaction of $(A,\de)$ on $D$.

Further, since $X$ is a corepresentation, we know that $X \in \M(A \ot
\te(\Ah))$. So, by the continuity of the coaction $\be$, we get,
\begin{align*}
[D \te(\Ah)] &= [(\om \ot \io)\bigl(X^* \be(z)_{21} X (\K(H) \ot
\te(\Ah)) \bigr) \mid \om \in \B(H)_* ] \\ &=  [(\om \ot \io)\bigl(X^*
\be(z)_{21} (\K(H) \ot \te(\Ah)) \bigr) \mid \om \in \B(H)_* ] \\ &=
 [(\om \ot \io)\bigl(X^* (\K(H) \ot B \te(\Ah)) \bigr) \mid \om \in \B(H)_* ] =
B \; .
\end{align*}
This ends the proof of the theorem.
\end{proof}

\begin{remark} \label{rem.strictcontinuity}
There is another way to characterize uniquely the \cst-algebra $D$. We
claim that $D$ is the unique \cst-subalgebra of $\M(B)^\be$ that
satisfies the following conditions.
\begin{itemize}
\item The map $\al : x \mapsto X^*(1 \ot x)X$ defines a continuous coaction
  of $(A,\de)$ on $D$.
\item The map $\al : \M(B)^\be \recht \M(\K(H) \ot D) : z \mapsto X^* (1 \ot
  z) X$ is well defined and continuous on the unit ball of $\M(B)^\be$
  if we equip $\M(B)^\be$ with the strict topology inherited from
  $\M(B)$ and $\M(\K(H) \ot D)$ with the strict topology.
\item $D \subset \M(B)^\be$ is non-degenerate in the sense that $B = [DB]$.
\end{itemize}
First observe that the \cst-algebra $D$ defined above satisfies these
conditions. Since $D$ obviously satisfies the first and third
condition, it
remains to prove the second condition.
But $\eta : B \recht \Ah \ltimesred D$ is a $^*$-isomorphism and
since the inclusion $\Ah \ltimesred D \recht \M(\K(H) \ot D)$ is
non-degenerate, we get a strictly continuous map $\eta : \M(B) \recht
\M(\K(H) \ot D)$. It suffices to restrict $\eta$ to $\M(B)^\be$.

We prove the uniqueness: suppose that $D_1$ and $D_2$ satisfy the
stated conditions. Then,
$$[D_1 D_2] = [(\om \ot \io)(\al(D_1)) \; D_2 \mid \om \in \B(H)_*]
\subset [(\om \ot \io)(\al(\M(B)^\be)) \; D_2 \mid \om \in \B(H)_*]
\subset D_2 \; .$$
On the other hand, let $(e_i)$ be a bounded approximate identity for the
\cst-algebra $D_1$. Since $D_1 \subset \M(B)^\be$ is non-degenerate,
we get that $(e_i)$ is a net in $\M(B)^\be$ that converges to $1$ in
the strict topology of $\M(B)$. Take $\om \in \B(H)_*$ such that
$\om(1)=1$. Then, $\bigl((\om \ot \io)\al(e_i)\bigr)_i$ is a net in
$\M(D_2)$ that converges strictly to $1$. It follows that $D_2 \subset
[D_1 D_2]$ because $(\om \ot \io)\al(e_i) \in D_1$ for all $i$. We
conclude that $D_2 = [D_1 D_2]$. By symmetry, we get $[D_2 D_1] =
D_1$. Taking the adjoint, we find that $D_1 = D_2$. This proves our claim.

When proving Theorem \ref{thm.existqhs} we have in a natural way that
$D \subset Q$ covariantly, where $Q$ is a von Neumann algebra on which $(M,\de)$
coacts. Moreover, we have that $\Ah \ltimesred D$ is a dense
subalgebra of $\Mh \ltimes Q$. We claim that this implies that $D$ is
dense in $Q$. Indeed, if we denote by $\alh : \Ah \ltimesred D \recht
\M((\Ah \ltimesred D) \ot \Ah)$ the dual coaction, we have
$$\al(D)= [(\io \ot \io \ot \om)\bigl( \Wh_{13} \alh(\Ah \ltimesred D)
\Wh^*_{13} \bigr) \mid \om \in \B(H)_*] \; .$$
Hence, the $\si$-weak closure of $\al(D)$ is equal to the $\si$-weak closure of
$$[(\io \ot \io \ot \om)\bigl( \Wh_{13} \alh(\Mh \ltimes Q)
\Wh^*_{13} \bigr) \mid \om \in \B(H)_*]$$
and so, equal to the $\si$-weak closure of $[\al((\om \ot \io)\al(Q)) \mid
\om \in \B(H)_*]$. Since $\al(Q) (\B(H) \ot 1)$ is $\si$-weakly dense
in $\B(H) \ot Q$, we conclude that $D$ is $\si$-weakly dense in $Q$.
\end{remark}

\section{Induction of coactions and Green imprimitivity} \label{sec.indcoac}

We fix a locally compact quantum group $(M,\de)$. {\it We suppose
  throughout this section that $(M,\de)$ is strongly regular.}
We fix a closed quantum subgroup $(M_1,\de_1)$. So, we have $\pih :
  \Mh_1 \recht \Mh$.

Suppose that $\eta : C \recht \M(A_1 \ot C)$ is a continuous coaction
of $(A_1,\de_1)$ on a \cst-algebra $C$. We want to define an induced
\cst-algebra $\Ind C$ with a continuous coaction $\Ind \eta$ of
$(A,\de)$ on $\Ind C$. Of course, when $C=\C$ with the trivial
coaction, we want to find again $D$ with the coaction $\de$ of
$(A,\de)$ by left translations on $D$.

We defined $D$ as a suitable \cst-subalgebra of $Q = M^\al$, where
$\al : M \recht M \ot M_1$ is the coaction of $(M_1,\de_1)$ on $M$ by
right translations. To define $C$, we have again at our disposal a \cst-algebra which
is too big and inside which we want to find $\Ind C$.

\begin{notation}
We denote
$$\Ctil = \{ X \in \M(\K(H) \ot C) \mid X \in (M' \ot 1)'
\quad\text{and}\quad (\al \ot \io)(X) = (\io \ot
\eta)(X) \} \; .$$
We equip $\Ctil$ with the strict topology inherited from $\M(\K(H) \ot C)$
and call this the strict topology of $\Ctil$.
\end{notation}

Remark that the expression $X \in \M(\K(H) \ot C) \cap (M' \ot 1)'$ is the necessarily awkward way of saying that $X \in M \ot \M(C)$ in some loose
sense. Observe that, when $C = \C$, we have $\Ctil = M^\al = Q$.

\begin{theorem} \label{thm.existsinduced}
There exists a unique \cst-subalgebra $\Ind C$ of $\Ctil$ that
satisfies the following conditions.
\begin{itemize}
\item $\de \ot \io : \Ind C \recht \M(A \ot \Ind C)$ defines a
  continuous coaction of $(A,\de)$ on $\Ind C$.
\item $\de \ot \io : \Ctil \recht \M(\K(H) \ot \Ind C)$ is well
  defined and strictly continuous on the unit ball of $\Ctil$.
\item $\Ind C \subset \Ctil$ is non-degenerate in the sense that $H
  \ot C = [(\Ind C) (H \ot C)]$.
\end{itemize}
We define $\Ind \eta := \de \ot \io$ and call it the induced coaction
of $\eta$.
\end{theorem}

The following theorem shows that our definition of $\Ind C$ is the
correct one.

\begin{theorem} \label{thm.cstimprim}
There exist canonical covariant Morita equivalences
$$\Ahu \ltimesfull \Ind C \Morita \Ahu_1 \ltimesfull C
\quad\text{and}\quad
\Ah \ltimesred \Ind C \Morita \Ah_1 \ltimesred C \; .$$
The covariance is understood with respect to the dual coactions on the
crossed products.
\end{theorem}

The rest of this section will consist in proving both theorems. We start by performing again the induction procedure as in Section
\ref{sec.induction}, but taking into account systematically a \cst-algebra $C$ on which is coacted by $(A_1,\de_1)$.

Fix a coaction $\eta : C \recht \M(A_1 \ot C)$ of $(A_1,\de_1)$ on a
\cst-algebra $C$.
Let $\E$ be a \cst-$B$-module and let $(X,\te)$ be a covariant pair
for $\eta$ consisting of a corepresentation $X \in \cL(A_1 \ot \E)$
and a representation $\te : C \recht \cL(\E)$.

Let $\F=\Ind \E$ be the induced \cst-$B$-module and let $Y = \Ind X$ be the
induced corepresentation of $(A,\de)$ on $\F$. We claim that there
exists a canonical strict $^*$-homomorphism $\tetil : \Ctil \recht
\cL(\F)$ which is covariant in the sense that
\begin{equation} \label{eq.piep}
(\io \ot \tetil)(\de \ot \io)(z) = Y^* (1 \ot \tetil(z)) Y
\quad\text{for all}\;\; z \in \tetil \; .
\end{equation}
In order to give a meaning to the previous equality, we have to be careful. We consider
$$\Ctil_1 = \{ x \in \M(\K(H) \ot \K(H) \ot C) \mid x \in (1 \ot M'
\ot 1)' \quad\text{and}\quad (\io \ot \al \ot
\io)(x) = (\io \ot \io \ot \eta)(x) \} \; .$$
The algebra $\Ctil_1$ plays the role of $\B(H) \ot \Ctil$. It is not
difficult to define $\io \ot \tetil$ as a strict $^*$-homomorphism
$\Ctil_1 \recht \cL(H \ot \F)$ (see Lemma \ref{lem.extendstrict} for
a related result) . On the other hand, we have $\de \ot
\io : \Ctil \recht \Ctil_1$. As a composition of both, the left hand
side of \eqref{eq.piep} makes sense.

In the induction procedure for corepresentations, an important role is
played by the \wst-imprimitivity bimodule $\cI$ defined in Definition
\ref{def.wstimprim}. We extend $\cI$ as follows.

\begin{notation}
We define
$$\cJ = \{ x \in \cL(H_1 \ot C,H \ot C) \mid (\pih' \ot \io)(V_1)_{12}
x_{13} V^*_{1,12} = (\io \ot \eta)(x) \} \; .$$
Observe that $\cI \ot 1 \subset \cJ$. We also define
$$P_1 = \{ x \in \M(\K(H_1) \ot C) \mid V_{1,12} x_{13}
V^*_{1,12} = (\io \ot \eta)(x) \}$$
and
$$P = \{ x \in \M(\K(H) \ot C) \mid V_{12} x_{13}
V^*_{12} = (\io \ot \eta)(x) \} \; .$$
Then, $\cJ$ is a $P$-$P_1$-bimodule, $\cJ^* \cJ \subset P_1$ and $\cJ
\cJ^* \subset P$. Observe also that $\Ctil = P \cap (M' \ot 1)'$.
\end{notation}

Recall that the induced \cst-$B$-module $\F$ is defined by
$H \ot \F \cong \cI \rot{\pil} (H \ot \E) \; ,$
where $\pil : \Mh_1 \recht \cL(H \ot \E)$ is the strict
$^*$-homomorphism defined by $(\io \ot \pil)(W_1) = (\io \ot
\pih)(W_1)_{12} X_{13}$.

We extend $\pil$ to $P_1$ as follows. Let $z \in P_1$. It is easy to check
that
$$V_{1,12} \; \bigl( X(\io \ot \te)(z) X^*\bigr)_{13} V^*_{1,12} =
\bigl( X(\io \ot \te)(z) X^*\bigr)_{13} \; .$$
So, $X(\io \ot \te)(z) X^* \in \cL(H_1 \ot \E) \cap (\Mh_1' \ot
1)'$. Hence, we expect to be able to apply $\pih \ot \io$ to $X(\io
\ot \te)(z) X^*$. Exactly as in the definition of $\pil$ we are a
little bit more careful and we define the strict $^*$-homomorphism
$$\pitill : P_1 \recht \cL(H \ot \E) : \pitill(z) (v \ot 1) \xi = (v \ot
1) \; X(\io \ot \te)(z) X^* \; \xi $$ for all $z \in P_1, \xi \in H_1 \ot \E$ and $v \in \B(H_1,H)$ intertwining $\pih$. So, we can define the
\cst-$B$-module $\cJ \rot{\pitill} (H \ot \E)$ and the following $B$-linear inclusion that preserves inner products.
$$\cI \rot{\pil} (H \ot \E) \hookrightarrow \cJ \rot{\pitill} (H \ot
\E) \; .$$
We claim that this inclusion is unitary. So, we have to prove that the
image of this inclusion is dense. For this it suffices to show that
$(\cI \ot 1)P_1$ is strictly dense in $\cJ$. But, $(\cI^* \ot 1)\cJ
\subset P_1$, which implies that $(\cI \ot 1) P_1$ contains $(\cI \cI^*
\ot 1)\cJ$. Now, $\cI \cI^*$ is weakly dense in $\Mh \ltimes Q$ and in
particular, $1$ can be approximated by elements in $\cI \cI^*$ and we
are done.

As a conclusion, we have
$$H \ot \F \cong \cI \rot{\pil} (H \ot \E) = \cJ \rot{\pitill} (H \ot
\E) \; .$$
This shows that we have a natural strict $^*$-homomorphism $P \recht
\cL(H \ot \F)$.

Although $P$ is not a von Neumann algebra, we have defined, in a
sense, a $B$-correspondence
$$\corresp{H \ot \F}{P}{\Mh} \; .$$
Finally, we have to turn this $B$-correspondence into a bicovariant
$B$-correspondence. In Section \ref{sec.induction}, we already defined
the product coaction $\al_{\Ftil}$ on $\Ftil = H \ot \F$. We want to show that this coaction
is covariant with respect to the representation $P \recht \cL(H \ot
\F)$ together with the right coaction $z \mapsto \Vh_{13} (z \ot 1)
\Vh_{13}^*$ of $(\Ah,\deh)$ on $P$.

For this, we have to observe that
$$z \mapsto \Vh_{13} (z \ot 1) (\io \ot \pih)(\Vh_1^*)_{13} \; , \quad
z \in \cJ$$
defines a right coaction of $(\Ah,\deh)$ on $\cJ$. Then, the inclusion
$\cI \ot 1 \subset \cJ$ is compatible with the coaction $\al_{\cI}$
defined on $\cI$ in Section \ref{sec.induction}. Hence, the coaction
$\al_{\Ftil}$ can be considered as a product coaction on $\cJ
\rot{\pitill} (H \ot \E)$. Then, it is clear that we get a bicovariant
$B$-correspondence
$$\bicorresp{H \ot \F}{P}{\Mh}{M'} \; .$$
Since $\Ctil \subset P$ and since $\Vh_{13} (z \ot 1)
\Vh_{13}^* = z \ot 1$ for all $z \in \Ctil$, we have found a strict
$^*$-homomorphism
$$\Ctil \recht \cL(\F)$$
which is covariant with respect to the induced corepresentation $Y$.

By now it should be clear that the following lemma can be proved
in the same way as Theorem \ref{thm.firstimprim}.

\begin{lemma} \label{lem.firstimprim}
Let $(A_1,\de_1)$ be a closed quantum subgroup of a strongly regular
locally compact quantum group $(A,\de)$. Let $\eta : C \recht \M(A_1
\ot C)$ be a continuous coaction of $(A_1,\de_1)$ on a \cst-algebra $C$.

A corepresentation $Y$ of $(A,\de)$ on a \cst-$B$-module $\F$ is
induced from a covariant pair $(X,\te)$ consisting of a
corepresentation $X$ of $(A_1,\de_1)$ and a representation of $C$ if
and only if there exists a strict $^*$-homomorphism $\Ctil \recht
\cL(F)$ which is covariant with respect to $Y$.
\end{lemma}

A final ingredient to prove Theorems
\ref{thm.existsinduced} and \ref{thm.cstimprim} is
Proposition \ref{prop.inducedcoaction}. Let $\be : \E \recht \M(\E \ot
\Ah)$ be a coaction of $(\Ah,\deh)$ on $\E$ which is compatible with a
continuous coaction on $B$ still denoted by $\be$. Suppose that $\be$
and the corepresentation $X \in \cL(A_1 \ot \E)$ are covariant in the
sense of \eqref{eq.nogcovar}. Then,
Proposition \ref{prop.inducedcoaction} yields an induced coaction
$\Ind \be$ on $\F$. Also the following lemma can be proved easily.

\begin{lemma} \label{lem.inductionofacoaction}
If $\te(C)$ is part of the fixed point algebra of $\cL(\E)$ with
respect to $\be$, then $\Ctil$ is part of the fixed
point algebra of $\cL(\F)$ with respect to $\Ind \be$.
\end{lemma}

We have then gathered enough material to prove Theorems
\ref{thm.existsinduced} and \ref{thm.cstimprim}.

\begin{proof}[Proof of Theorem \ref{thm.existsinduced}]
Consider the \cst-$\Ah_1 \ltimesred C$-module $\Ah_1 \ltimesred C$
together with the reduced covariant pair $(X,\te)$ consisting of the
corepresentation $X = W_{1,12} \in \cL(A_1 \ot (\Ah_1 \ltimesred C))$
and the representation $\te = \eta : C \recht \M(\Ah_1 \ltimesred C)$.

Define $\F = \Ind (\Ah_1 \ltimesred C)$, together with the induced
corepresentation $Y \in \cL(A \ot \F)$ and the strict  \linebreak
$^*$-homomorphism $\Ctil \recht \cL(\F)$.

On the crossed product $\Ah_1 \ltimesred C$, we have the dual coaction
that we push into $(\Ah,\deh)$ using $\pih$:
$$\etah : \Ah_1 \ltimesred C \recht \M\bigl((\Ah_1 \ltimesred C) \ot
\Ah) \; .$$
By Lemma \ref{lem.inductionofacoaction}, we find a coaction $\ga$ on
$\F$ such that $(\io \ot \ga)(Y) = W_{13} Y_{12}$ and such that
$\Ctil$ is part of the fixed point algebra $\cL(\F)^\ga$.

By Theorem \ref{thm.quantumlandstad} and Remark
\ref{rem.strictcontinuity}, we find a \cst-algebra $\Ind C
\subset \cL(\F)^\ga$ such that
\begin{itemize}
\item $z \mapsto Y^*(1 \ot z)Y$ defines a continuous coaction of
  $(A,\de)$ on $\Ind C$;
\item $\cL(\F)^\ga \recht \M(\K(H) \ot \Ind C) : z \mapsto Y^*(1 \ot
  z)Y$ is well defined and strictly continuous;
\item $\K(\F) = [(\Ind C) \; (\om \ot \io)(Y) \mid \om \in \B(H)_*]$.
\end{itemize}

Consider then the regular representation $\Ah_1 \ltimesred C \recht
\cL(H_1 \ot C)$. It is straightforward to check that $\Ind
(H_1 \ot C) = H \ot C$. Since
$$(\Ind (\Ah_1 \ltimesred C)) \rot{\Ah_1 \ltimesredsmall C} (H_1 \ot C) =
\Ind(H_1 \ot C) = H \ot C\; ,$$
we find a faithful $^*$-homomorphism $\mu : \cL(\F) \recht \cL(H \ot
C)$.

One verifies that $\mu(\cL(\F)) \subset P$. Further, $\mu$ intertwines
the coaction $\ga$ on $\cL(\F)$ with the coaction $z \mapsto \Vh_{13}
(z \ot 1) \Vh_{13}^*$ of $(\Ah,\deh)$ on $P$. Hence, $\mu\bigl(
\cL(\F)^\ga \bigr) \subset \Ctil$. So, we have shown that $\cL(\F)^\ga
\cong \Ctil$ in a canonical way and that the isomorphism preserves the
strict topology on bounded subsets.

As such, we have defined
$\Ind C$ in a canonical way as a subalgebra of $\Ctil$ such that $\Ind
C$ satisfies the required conditions.
The uniqueness of $\Ind C$ is proved in exactly the same way as the
uniqueness statement in Remark \ref{rem.strictcontinuity}.
\end{proof}

We finally prove Theorem \ref{thm.cstimprim}.

\begin{proof}[Proof of Theorem \ref{thm.cstimprim}]
In the proof of Theorem \ref{thm.existsinduced} the reduced Morita
equivalence has already been shown. We shall prove the full Morita equivalence.

We can almost copy the proof of Theorem \ref{thm.morita}. The only
point to show is the following: let $Y \in \cL(A \ot \F)$ be a
corepresentation of $(A,\de)$ on a \cst-$B$-module $\F$ which is
covariant with respect to a representation $\rho : \Ind C \recht
\cL(\F)$, then $\rho$ extends uniquely to a strict $^*$-homomorphism
$\Ctil \recht \cL(\F)$ which remains covariant. Once we have proved
this claim, we can use Lemma \ref{lem.firstimprim}.

The uniqueness of the extension is obvious from the requirement of
covariance and from the existence of the map $\de \ot \io : \Ctil \recht \M(\K(H) \ot \Ind
C)$.

To prove the existence of the extension of $\rho$, we define
$$\mu : \Ctil \recht \cL(H \ot \F) : \mu(z) = Y (\io \ot \rho)(\de \ot
\io)(z) Y^* \; .$$
Observe that $\mu(z) = 1 \ot \rho(z)$ for all $z \in \Ind C$. Let $z
\in \Ctil$. Clearly, $\mu(z) \in (M' \ot 1)'$. Moreover,
$$(\de \ot \io)\mu(z) = Y_{13} Y_{23} \; (\io \ot \io \ot \rho)(\io
\ot \de \ot \io)(\de \ot \io)(z) \; Y^*_{23} Y^*_{13} \; .$$
Since $(\de \ot \io)(z) \in \M(\K(H) \ot \Ind C)$, the covariance of
$\rho$ implies that
$$Y_{23} \; (\io \ot \io \ot \rho)(\io
\ot \de \ot \io)(\de \ot \io)(z) \; Y^*_{23} = \bigl((\io \ot
\rho)(\de \ot \io)(z)\bigr)_{13} \; .$$
It follows that $(\de \ot \io)\mu(z) = \mu(z)_{13}$, which implies
that $\mu(z) = 1 \ot \rho_1(z)$ for a well defined strict \linebreak
$^*$-homomorphism $\rho_1 : \Ctil \recht \cL(\F)$. Then, $\rho_1$ is
the required extension.
\end{proof}

\begin{remark} \label{rem.reduce-not}
In order to define $\Ind C$, we did not assume that the coaction of $(A_1,\de_1)$ on $C$ is reduced. One checks easily that the induction of
$(C,\eta_1)$ coincides with the induction of the reduction of $(C,\eta_1)$. This is not very surprising, since our induced coaction is reduced by
definition. Nevertheless, we prefer not to impose that $\eta_1$ is reduced, to get neater formulations in the next section.
\end{remark}

\section{Induction and restriction} \label{sec.indres}

We still fix a strongly regular locally compact quantum group
$(M,\de)$ with a closed quantum subgroup $(M_1,\de_1)$. Suppose that
$\eta : C \recht \M(A \ot C)$ is a continuous coaction of $(A,\de)$ on
the \cst-algebra $C$. We shall deal with the following problem: what
happens if we first restrict the coaction $\eta$ to a coaction $\eta_1 : C
\recht \M(A_1 \ot C)$ and then induce $\eta_1$?

In the classical situation, where $G_1$ is a closed subgroup of a
locally compact group $G$ and where $G$ is acting continuously on $C$,
we know that the induction of the restriction will be $C_0(G/G_1) \ot
C$ with the diagonal action. There is an obvious reason for which we
cannot expect exactly the same result: if $M$ is non-commutative, the
notion of a diagonal coaction does not make sense in general. We shall
see however in Proposition \ref{thm.inductionrestriction} that there is a
natural description of the induction of the restriction of a coaction,
where the tensor product with diagonal action is replaced by a twisted
product and diagonal coaction.

Of course we should first discuss the construction of the restriction
of $\eta$ to $\eta_1$. Recall that the morphism
$(M,\de) \overset{\pi}{\longrightarrow} (M_1,\de_1)$ comes as a
non-degenerate $^*$-homomorphism $\pi : \Au \recht \M(\Au_1)$. Hence,
the obvious formula $\eta_1 = (\pi \ot \io)\eta$ does not make sense
immediately.

\begin{definition}
Let $\eta : C \recht \M(A \ot C)$ be a continuous coaction which
admits a lift to a coaction $C \recht \M(\Au \ot \C)$. Then the restriction $\eta_1 : C \recht
\M(A_1 \ot C)$ is the unique coaction satisfying
$$(\io \ot \eta_1)\eta = (\al \ot \io)\eta \; ,$$
where $\al : A \recht \M(A \ot A_1)$ is the right coaction of
$(A_1,\de_1)$ on $(A,\de)$. The restriction $\eta_1$ is continuous.
\end{definition}
Observe that $\eta$ admits a unique lift to a coaction $C \recht \M(\Au \ot C)$ when either $\eta$ is reduced or maximal, but that an arbitrary
continuous coaction does not necessarily admit such a lift.

\begin{theorem} \label{thm.inductionrestriction}
Let $\eta : C \recht \M(A \ot C)$ be a continuous coaction of
$(A,\de)$ on the \cst-algebra $C$, admitting a lift to a coaction $C
\recht \M(\Au \ot C)$. Restrict $\eta$ to a coaction
$\eta_1 : C \recht \M(A_1 \ot C)$. The induced \cst-algebra for
$\eta_1$ is given by
$$\Ind(C,\eta_1) \cong [(D \ot 1) \eta(C)]$$
with induced coaction $\de \ot \io$.

Moreover, $D \recht \M(D \ot \Ah) : x \mapsto W(x \ot 1)W^*$ defines a
continuous right coaction of $(\Ah,\dehop)$ on $D$. We call this the
\emph{adjoint coaction}. If $\eta$ is a reduced coaction, we have
$$\Ind(C,\eta_1) \cong [ \; (W(D \ot 1)W^*)_{12} \; \eta(C)_{23} \; ] \;
.$$
\end{theorem}

Remark that in the case where $M$ is commutative, the adjoint coaction
is trivial and we find that $\Ind(C,\eta_1) \cong D \ot C$.

In order to prove Theorem \ref{thm.inductionrestriction}, we have to
do some preliminary work.
We have shown in Theorem \ref{thm.morita} that $\Ah \ltimesred D
\Morita \Ah_1$. In fact, the Morita equivalence is given by $\F = \Ind
\Ah_1$, which means that
$$H \ot \F \cong \cI \rot{\pil} (H \ot \Ah_1) \; .$$
We have seen that $\cI \rot{\pil} (H \ot \Ah_1)$ is a bicovariant
$\Ah_1$-correspondence:
$$\bicorresp{\cI \rot{\pil} (H \ot \Ah_1)}{\Mh \ltimes Q}{\Mh}{M'} \;
.$$
Consider now the \wst-$\Mh_1$-module $H \ot \cI$. Then, equipped with
the obvious representations of $\Mh \ltimes Q$, $\Mh' \ot
1$, $M' \ot 1$, we get an inclusion
$$\cI \rot{\pil} (H \ot \Ah_1) \hookrightarrow H \ot \cI : v
\rot{\pil} \xi \mapsto \Vh_{21} (1 \ot v) (\pih \ot \io)(\Vh_{1,21}^*)
\xi$$
intertwining these three representations. It follows that there is a
closed subspace $\cI_0 \subset \cI$ such that (applying a flip map)
\begin{itemize}
\item $\cI_0 \ot H = [\Vh (\cI \ot 1) (\io \ot \pih)(\Vh_1^*) (\Ah_1
  \ot H)]$;
\item $[\cI_0^* \cI_0] = \Ah_1$;
\item $[\cI_0 \cI_0^*] = [\Ah D]$;
\item $v \mapsto \Vh (v \ot 1)(\io \ot \pih)(\Vh_1^*)$ defines a
  continuous coaction of $(\Ah,\deh)$ on $\cI_0$.
\end{itemize}
We should consider $\cI_0$ as a \emph{concrete realization of the $\Ah
  \ltimesred D$-$\Ah_1$ imprimitivity bimodule}.

Suppose now that $\eta_1 : C \recht \M(A_1 \ot C)$ is a continuous
coaction of $(A_1,\de_1)$ on the \cst-algebra $C$. We can give an
analogous concrete $\Ah \ltimesred \Ind C$ - $\Ah_1 \ltimesred C$
imprimitivity bimodule. From the proof of Theorem
\ref{thm.existsinduced}, we know that such an imprimitivity bimodule
is given by $\Ind (\Ah_1 \ltimesred C)$. In exactly the same way as
above, we then find that
$$\Ind(\Ah_1 \ltimesred C) \cong [(\cI_0 \ot 1)\eta_1(C)] \subset
\M(\K(H_1,H) \ot C)$$
as a covariant \cst-bimodule.

\begin{proof}[Proof of Theorem \ref{thm.inductionrestriction}]
We first prove the claim on the adjoint coaction. On the
\wst-$\Mh_1$-module $\cI$ we consider the coaction $\cI \recht \cI \ot \Mh :
v \mapsto W(v \ot 1)(\io \ot \pih)(W_1)^*$ of $(\Mh,\dehop)$
on the right, which is compatible with
the coaction $(\io \ot \pih)\dehopo : \Mh_1 \recht \Mh_1 \ot \Mh$ of
$(\Mh,\dehop)$ on $\Mh_1$. We also consider the coaction $\io \ot (\io
\ot \pih)\dehopo$ of $(\Ah,\dehop)$ on $H \ot \Ah_1$. From
Proposition \ref{prop.productofcoactions} we get a product coaction of
$(\Ah,\dehop)$ on $\cI \rot{\pil} (H \ot \Ah_1)$. Since it is clear
that this product coaction leaves invariant the representations of
$M'$ and $\Mh'$, we get a coaction
$$\ga : \Ind \Ah_1 \recht \M((\Ind \Ah_1) \ot \Ah)$$
of $(\Ah,\dehop)$ on the imprimitivity bimodule $\Ind \Ah_1$, which is
compatible with the coaction $(\io \ot \pih)\dehopo$ of $(\Ah,\dehop)$
on $\Ah_1$. If we follow the isomorphism $\Ind \Ah_1 \cong \cI_0$, the
coaction $\ga$ is given by $\ga : \cI_0 \mapsto \M(\cI_0 \ot \Ah) : v
\mapsto W(v \ot 1)(\io \ot \pih)(W_1^*)$.

In particular, we get a coaction $\ga :\Ah
\ltimesred D \recht \M((\Ah
\ltimesred D) \ot \Ah)$.
By construction, the coaction $\ga$ commutes with the
dual coaction on $\Ah \ltimesred D$ and is given by $\dehop$ on
$\Ah$. This implies (see the proof of Theorem
\ref{thm.quantumlandstad}) that the restriction
of $\ga$ to $D$ is a continuous coaction of $(\Ah,\dehop)$ on $D$. By
construction, this coaction is exactly given by $x \mapsto W(x \ot
1)W^*$.

Above we have concretely realized the $\Ah \ltimesred \Ind C$ - $\Ah_1 \ltimesred C$
imprimitivity bimodule as $[(\cI_0 \ot 1)\eta_1(C)]$.
We claim that $[(\cI_0 \ot 1) \eta_1(C)] = [\eta(C) (\cI_0 \ot
1)]$. To prove this claim, we write $\etatil(x) = (J\Jh \ot 1) \eta(x) (\Jh J
\ot 1)$ and observe that
$$(\io \ot \eta_1) \etatil(x) = (\pih \ot \io)(\Wh_1^*)_{12}
\etatil(x)_{13} (\pih \ot \io)(\Wh_1)_{12} \quad\text{and}\quad
(\io \ot \eta)\etatil(x) = \Wh^*_{12} \etatil(x)_{13} \Wh_{12} \; .$$
But then, by continuity of the coaction $\ga : \cI_0 \recht \M(\cI_0
\ot \Ah)$ defined above, we get
\begin{align*}
[\K \ot (\cI_0 \ot 1) \eta_1(C)] &= [(\K \ot \cI_0 \ot 1) (\io \ot
\eta_1) \etatil(C)] \\ &= [((1 \ot \K)\ga(\cI_0) \ot 1) (\io \ot
\pih)(W_1)_{12} \etatil(C)_{23} (\io \ot \pih)(W_1^*)_{12}]_{213} \\
&= [(1 \ot \K \ot 1) W_{12} (\cI_0 \ot \etatil(C)) (\io \ot
\pih)(W_1^*)_{12}]_{213}
\\ &= [(\K \ot 1 \ot 1) (\io \ot \eta)\etatil(C) \ga(\cI_0)_{21} ] \\
&= [(\K \ot \eta(C)) \ga(\cI_0)_{21}] = \K \ot [\eta(C) (\cI_0 \ot 1)]
\; .
\end{align*}
This proves our claim. Since $\Ah \ltimesred \Ind C \cong \K\bigl(
\Ind (\Ah_1 \ltimesred C) \bigr)$, we conclude that
$$\Ah \ltimesred \Ind C \cong [(\cI_0 \ot 1) \eta_1(C) (\cI_0^* \ot 1)
\eta(C)] = [(\cI_0 \cI_0^* \ot 1) \eta(C)] = [(\Ah D \ot 1) \eta(C)]
$$
in such a way that the dual coaction on $\Ah \ltimesred \Ind C$
corresponds with the coaction $z \mapsto \Vh_{13} (z \ot 1)
\Vh_{13}^*$ of $(\Ah,\deh)$ on $[(\Ah D \ot 1) \eta(C)]$. It follows
from Theorem \ref{thm.quantumlandstad} that $\Ind C \cong [(D \ot
1)\eta(C)]$ in such a way that the induced coaction $\Ind \eta_1$
corresponds to $\de \ot \io$.

We finally observe that, when $\eta$ is a reduced coaction
$$[(W(D \ot 1)W^*)_{12} \eta(C)_{23}] = W_{12} (\io \ot \eta)\bigl(
[(D \ot 1) \eta(C)] \bigr) W_{12}^* \cong [(D \ot 1)\eta(C)] \; .$$
This ends the proof of the theorem.
\end{proof}

We give a slightly more natural explanation why $[(D \ot 1)\eta(C)]$ is indeed a \cst-algebra. It is an example of a kind of \emph{twisted product}
generalizing reduced crossed products.

\begin{proposition} \label{prop.twisted}
Let $(A,\de)$ be a locally compact quantum group and let $C,D$ be \cst-algebras. Suppose that $\ga : D \recht \M(\Ah \ot D)$ is a continuous coaction
of $(\Ah,\deh)$ on $D$ and that $\al : C \recht \M(A \ot C)$ is a continuous coaction of $(A,\de)$ on $C$. Then,
$$[\ga(D)_{12} \; \al(C)_{13}] \subset \M(\K(H) \ot D \ot C)$$
is a \cst-algebra.
\end{proposition}

This proposition indeed generalizes the notion of reduced crossed
product. If $\al : C \recht \M(A \ot C)$ is a continuous coaction, we
take $D = \Ah$ with the coaction $\ga = \deh$. We observe that, using
the notation $\si$ for the flip,
$$(\si \ot \io)[\deh(\Ah)_{12} \; \al(C)_{13}] = W_{12} [(\Ah \ot 1
\ot 1) (\io \ot \al)\al(C)] W_{12}^* = (WV)_{12}  (\Ah
\ltimesred C)_{13} (WV)_{12}^* \cong \Ah \ltimesred C \; .$$

To obtain $\Ind (C,\eta_1)$ as an example of such a twisted product, we
take the adjoint coaction on the quantum homogeneous space $D$ and the
given coaction $\eta : C \recht \M(A \ot C)$ on $C$.

\begin{proof}
By the continuity of $\ga$, we get
$$\ga(D) = [(\om \ot \io \ot \io)(\deh \ot \io)\ga(D)] = [(\om \ot \io
\ot \io)(\Vh_{12} \ga(D)_{13} \Vh^*_{12})] = [(\om \ot \io \ot
\io)(W^*_{12} \gatil(D)_{13} W_{12})] \; ,$$
where $\gatil(x) = (J\Jh \ot 1)\ga(x)(\Jh J \ot 1)$. Using the
continuity of $\al$, it follows that
\begin{align*}
[\ga(D)_{12} \; \al(C)_{13}] &= [(\om \ot \io \ot \io \ot
\io)(W^*_{12} \gatil(D)_{13} W_{12} \; (\K(H) \ot \al(C))_{124} )] \\ &=
[(\om \ot \io \ot \io \ot
\io)(W^*_{12} \gatil(D)_{13} W_{12} \; (\de \ot \io)\al(C)_{124})] \\
&= [(\om \ot \io \ot \io \ot
\io)(W^*_{12} \gatil(D)_{13} \; \al(C)_{24} W_{12})] \\
&= [(\om \ot \io \ot \io \ot
\io)( (\de \ot \io)\al(C)_{124} \; W^*_{12} \gatil(D)_{13} W_{12})] \\
&= [(\om \ot \io \ot \io \ot
\io)( \al(C)_{24} \; W^*_{12} \gatil(D)_{13} W_{12})] = [\al(C)_{13}
\; \ga(D)_{12}] \; .
\end{align*}
It is then clear that $[\ga(D)_{12} \; \al(C)_{13}]$ is a \cst-algebra.
\end{proof}

As an application of these twisted products, we define \emph{reduced crossed
  products by homogeneous spaces}.

\begin{proposition} \label{prop.crossedcoideal}
Let $(A,\de)$ be a strongly regular locally compact quantum group with
a closed quantum subgroup $(A_1,\de_1)$. Denote the quantum
homogeneous space by $D$.

Let $\mu : C \recht \M(\Ah \ot C)$ be continuous coaction of
$(\Ah,\deh)$ on a \cst-algebra $C$. Then,
$$[(D \ot 1) \mu(C)]$$
is a \cst-algebra that we denote by $D \ltimesred C$ and that we call
the \emph{reduced crossed product} of $C$ by $D$.
\end{proposition}
\begin{proof}
It suffices to observe that $\de : D \recht \M(A \ot D)$ is a
continuous coaction of $(A,\de)$ on $D$. By Proposition
\ref{prop.twisted} we get that $[\de(D)_{12} \mu(C)_{13}]$ is a
\cst-algebra. But,
$$[\de(D)_{12} \mu(C)_{13}] = W^*_{12} (\si \ot \io)[(D \ot 1 \ot 1)
(\io \ot \mu)\mu(C)] W_{12} = W^*_{12} \Vh_{21} (1 \ot [(D \ot
1)\mu(C)]) \Vh_{21}^* W_{12} \; .$$
It follows that $[(D \ot 1)\mu(C)]$ is a \cst-algebra.
\end{proof}

Finally it is possible to define as well \emph{full crossed products
  by homogeneous spaces}. Nevertheless we should remark that our
  definition does not fully generalize full crossed products by
  quantum groups (see Remark \ref{rem.opgepast} below).

\begin{definition} \label{def.fullqhs}
Let $(A,\de)$ be a strongly regular locally compact quantum group with
a closed quantum subgroup $(A_1,\de_1)$. Denote the quantum
homogeneous space by $D$. Let $\mu : C \recht \M(\Ah \ot C)$ be continuous coaction of
$(\Ah,\deh)$ on a \cst-algebra $C$.

A pair $(\pi_D,\pi_C)$ of non-degenerate representations of $D$ and
$C$ on a Hilbert space $K$ is said to be \emph{covariant} if the
\cst-algebras
$$(\io \ot \pi_D)(V(1 \ot D)V^*) \quad\text{and}\quad (\io \ot
\pi_C)\mu(C)$$
commute as \cst-algebras on $H \ot K$.

If $(\pi_D,\pi_C)$ is such a covariant representation, then $[\pi_D(D)
\pi_C(C)]$ is a \cst-algebra. The \cst-algebra generated by a
universal covariant representation is denoted by $D \ltimesfull C$ and
called the \emph{full crossed product of $C$ by $D$}.
\end{definition}

Observe that the expression $(\io \ot \pi_D)(V(1 \ot D)V^*)$ makes
sense because the adjoint coaction is a well defined continuous
coaction on $D$.

\begin{remark} \label{rem.opgepast}
Suppose that, in the setting of the previous definition, $A_1=\C$, the
one-point subgroup. Then, of course, $D = A$. One should expect that
$D \ltimesfull C$ coincides with the usual full crossed product. This
is not necessarily the case: the representation $\pi_D$ is a
representation of $D=A$ and not of the full \cst-algebra $\Au$.

It is possible to define a \emph{universal homogeneous space} and to
use this one, rather than $D$, to define full crossed products. We do
not go into this: in the situation where we use the full
crossed product by $D$ (Theorem \ref{thm.mansfield}), the coaction $\mu$ is a
dual coaction and hence, there is no need to take a universal variant
of $D$.
\end{remark}

\section{Characterization of induced coactions}

It is well known that a continuous action of a l.c.\ group $G$ on a \cst-algebra $B$ is induced from an action of a closed subgroup $G_1$ if and only
if there exists a $G$-equivariant embedding of $C_0(G/G_1)$ into the center of $\M(B)$. We prove a similar result in this section. Nevertheless, we
cannot hope for an identical characterization since it would require the commutativity of the quantum homogeneous space $D$.

We fix a strongly regular l.c.\ quantum group $(M,\de)$ together with a closed quantum subgroup $(M_1,\de_1)$ given by the morphism $\pih : \Mh_1
\recht \Mh$. We denote as above the quantum homogeneous space by $D$, which is a \cst-subalgebra of the measured quantum homogeneous space $Q =
M^\al$, where $\al$ is the right coaction of $(M_1,\de_1)$ on $M$ by right translations.

\begin{lemma} \label{lem.antiauto}
There is a canonical anti-automorphism $\gamma$ of $\Ah \ltimesred D$ given by $\gamma(z) = Jz^*J$ when $\Ah \ltimesred D$ is represented on $H$.
\end{lemma}
\begin{proof}
Consider the concrete $\Ah \ltimesred D$ - $\Ah_1$ imprimitivity bimodule $\cI_0 \subset \B(H_1,H)$ constructed in the previous section. We claim
that $\cI_0 = J \cI_0 J_1$. It is clear that the lemma follows from this claim. First of all, $J \cI J_1 = \cI$, since $J \pih(x) J = \Rh(\pih(x^*))
= \pih(\Rh_1(x^*)) = \pih(J_1 x J_1)$ for all $x \in \Mh_1$. Hence, it follows that
\begin{align*}
J \cI_0 J_1 \ot H &= [W (J\cI J_1 \ot 1) (\io \ot \pih)(W_1^*) (\Ah_1 \ot H)] = [W (\cI \ot 1) (\io \ot \pih)(W_1^*) (\Ah_1 \ot H)] \\ &\supset [W
(\cI_0 \ot 1) (\io \ot \pih)(W_1^*) (\Ah_1 \ot H)] \; .
\end{align*}
In the proof of Theorem \ref{thm.inductionrestriction}, we have shown that the coaction $v \mapsto W(v \ot 1)(\io \ot \pih)(W_1^*)$ is continuous. It
then follows that $J \cI_0 J_1 \ot H \supset \cI_0 \ot H$, which proves our claim.
\end{proof}

\begin{theorem}
Let $\eta : B \recht \M(A \ot B)$ be a continuous reduced coaction of $(A,\de)$ on a \cst-algebra $B$. Then $\eta$ is induced from a coaction of
$(A_1,\de_1)$ if and only if there exists a non-degenerate $^*$-homomorphism $\te : D \recht \M(B)$ which is covariant and satisfies the condition
$$\tetil(\gamma(D)) \quad\text{commutes with}\quad B \quad\text{in}\quad \Ah \ltimesred B \; ,$$
where $\tetil : \Ah \ltimesred D \recht \M(\Ah \ltimesred B)$ denotes the extension of $\te$ and $\gamma$ is the anti-automorphism defined in Lemma
\ref{lem.antiauto}.
\end{theorem}

Observe that, in the case where $M=L^\infty(G)$ is abelian, we have $\gamma(D) = D$ and the above condition exactly says that $\te(D)$ is in the
center of $\M(B)$.

\begin{proof}
Consider again the concrete $\Ah \ltimesred D$ - $\Ah_1$ imprimitivity bimodule $\cI_0 \subset \B(H_1,H)$ constructed in the previous section. If
$\eta_1$ is a continuous coaction of $(A_1,\de_1)$ on $C$, we can realize $\Ind C \subset \M(A \ot C)$ such that
$$\Ah \ltimesred \Ind C = [(\Ah \ot 1) \Ind C] = [(\cI_0 \ot 1) \eta_1(C) (\cI_0^* \ot 1)] \; .$$
It then follows easily that we can embed $\te : D \recht \M(\Ind C) : \te(x) = x \ot 1$ such that $\tetil(\gamma(D)) = JDJ \ot 1 \subset M' \ot 1
\subset (\Ind C)'$.

Suppose conversely that we have a continuous coaction of $(A,\de)$ on $B$ and a $^*$-homomorphism $\te : D \recht \M(B)$ with the properties stated
in the theorem. Define $\E = \cI_0^* \rot{\tetil} (\Ah \ltimesred B)$ and denote by $\rho$ the coaction of $(\Ah,\deh)$ on $\E$ on the right, which
is the internal tensor product of the coaction $\cI_0^* \recht \M(\cI_0^* \ot \Ah) : v \recht (\io \ot \pih)(\Vh_1) (v \ot 1) \Vh^*$ and the dual
coaction on $\Ah \ltimesred B$. We still write $\rho$ for the corresponding coaction on $\K(\E)$. Observe that we have a representation $\Ah_1 \recht
\cL(\E)$ which is covariant with respect to the coaction $(\io \ot \pih_1)\deh_1$ of $(\Ah,\deh)$ on $\Ah_1$ and the coaction $\rho$ on $\cL(\E)$.

We write $F:=\K(\E)$. We claim that $\rho :F \recht \M(F \ot \Ah)$ is of the form $(\io \ot \pih)\rho_1$ for a continuous coaction $\rho_1$ of
$(\Ah_1,\deh_1)$ on $F$. It then follows from Theorem \ref{thm.quantumlandstad} that $F = \Ah_1 \ltimesred C$ for some continuous coaction $\eta_1$
of $(A_1,\de_1)$ on $C$. By construction we have $(B,\eta) \cong (\Ind C,\Ind \eta_1)$.

To prove our claim, take $x \in D$, $v \in \cI_0$ and $b \in B$. Then
\begin{align*}
(1 \ot \Jh x \Jh) \rho(v^* \rot{\tetil} b) & = (1 \ot \Jh x \Jh) \bigl( (\io \ot \pih)(\Vh_1) (v^* \ot 1) \Vh^* \rot{\tetil \ot \io} (b \ot 1) \bigr)
\\ &= \bigl((\io \ot \pih)(\Vh_1)(v^* \ot 1)\Vh^* \; (J \ot \Jh) \deop(x) (J \ot \Jh) \bigr) \rot{\tetil \ot \io} (b \ot 1)
\\ &= \bigl( (\io \ot \pih)(\Vh_1)(v^* \ot 1)\Vh^* \;  (\gamma \ot R)\deop(x^*) \bigr) \rot{\tetil \ot \io} (b \ot 1)
\\ &= (\io \ot \pih)(\Vh_1)(v^* \ot 1)\Vh^* \rot{\tetil \ot \io} (\tetil\gamma \ot R)\deop(x^*) (b \ot 1) = \rho(v^* \rot{\tetil} b) \; (\tetil\ga \ot
R)\deop(x^*) \; .
\end{align*}
It follows from the calculation that $1 \ot \Jh D \Jh$ and $\rho(F)$ commute. Since $M \cap (\Jh D \Jh)' = \pih(\Mh_1)$, it follows that $\rho(F)$
has its second leg in $\pih(\Mh_1)$. So, we can define $\rho_1 : F \recht \M(F \ot \K(H_1)) \cap (1 \ot \Mh_1')'$ such that $(1 \ot w) \rho_1(x) =
\rho(x)(1 \ot w)$ whenever $x \in F$ and $w \in \B(H_1,H)$ with $w y = \pih(y) w$ for all $y \in \Mh_1$. It follows that $[\rho_1(F) (1 \ot H_1)] = F
\ot H_1$. Since $\rho$ is a continuous coaction, we find that
\begin{align*}
[\rho_1(F)(1 \ot \Ah_1) \ot H_1] &= [(\rho_1 \ot \io)\rho_1(F) (1 \ot \Ah_1 \ot H_1)] = [\Wh_{1,23}^* \rho_1(F)_{13} \Wh_{1,23} (1 \ot \Ah_1 \ot
H_1)]
\\ &= [\Wh_{1,23}^* \rho_1(F)_{13}(1 \ot \Ah_1 \ot H_1)] = [\Wh_{1,23}^* (F \ot \Ah_1 \ot H_1)] = F \ot \Ah_1 \ot H_1 \; .
\end{align*}
This proves that $\rho_1$ is a continuous coaction of $(\Ah_1,\deh_1)$ on $F$ and hence, proves our claim and the theorem.
\end{proof}

\section{Green-Rieffel-Mansfield imprimitivity} \label{sec.mansfield}

We still fix a strongly regular locally compact quantum group
$(M,\de)$ with a closed quantum subgroup $(M_1,\de_1)$. So, we have
$\pih : \Mh_1 \recht \Mh$.

Suppose that $\eta : C \recht \M(JAJ \ot C)$ is a continuous coaction of
$(JAJ,\de')$ on a \cst-algebra $C$, which admits a lift to the
universal level.
(It is not crucial to take a
left coaction of $JAJ$ rather then a right coaction of $(A,\de)$. The
only convenience is that the crossed product is now $\Ah \ltimes C$ with
the dual action being a left coaction of $(\Ah,\deh)$ such that the
second crossed product is $A \ltimes \Ah \ltimes C$.) The
comultiplication $\de'$ is defined by
$$\de'(JxJ) = (J \ot J) \de(x) (J \ot J)$$
for all $x \in A$.

We then have a restricted coaction $\eta_1 : C \recht \M(J_1 A_1 J_1
\ot C)$.

\begin{theorem} \label{thm.mansfield}
If $\eta$ is a reduced coaction, there is a canonical Morita
equivalence
$$D \ltimesred \Ah \ltimesred C \Morita \Ah_1 \ltimesred C \; .$$
If $\eta$ is a maximal coaction, there is a canonical Morita
equivalence
$$D \ltimesfull \Ahu \ltimesfull C \Morita \Ahu_1 \ltimesfull C \; .$$
\end{theorem}

The conditions for $\eta$ being reduced or maximal are very natural. Indeed, if $(M_1,\de_1)$ is the one-point subgroup of $(M,\de)$, a natural
Morita equivalence between the second crossed product $A \ltimesred \Ah \ltimesred C$ and $C$ exists exactly when $\eta$ is reduced. On the other
hand, a natural Morita equivalence between the second crossed product $\Au \ltimesfull \Ahu \ltimesfull C$ and $C$ exists exactly when $\eta$ is
maximal.

The reduced and full crossed products by the quantum homogeneous space $D$ have
been defined in Proposition \ref{prop.crossedcoideal} and Definition \ref{def.fullqhs}.

Observe that $D \ltimesfull \Ahu \ltimesfull C$ is the universal \cst-algebra defined by covariant triples $(\rho,Y,\te)$ consisting of commuting
representations $\rho$ of $D$ and $\te$ of $C$, and a corepresentation $Y \in \M(A \ot \K(K))$ satisfying the covariance relations
\begin{equation}\label{eq.againcov}
(\io \ot \rho)\de(x) = Y^* (1 \ot \rho(x)) Y \quad\text{for all}\;\;
x \in D \quad\text{and}\quad (\io \ot \te)\etatil(y) = Y (1 \ot
\te(y))Y^* \quad\text{for all}\;\; y \in C \; .
\end{equation}
Here we used the following notation
$$\etatil : C \recht \M(A \ot C) : \etatil(y) = (J\Jh \ot
1)\eta(y)(\Jh J \ot 1) \; .$$
For later use, we also introduce the notation
$$\etatil_1 : C \recht \M(A_1 \ot C) : \etatil_1(y) = (J_1\Jh_1 \ot
1)\eta_1(y)(\Jh_1 J_1 \ot 1) \; .$$
Then, $\etatil$ is a coaction of $(A,\deop)$ on $C$ and $\etatil_1$ is
its restriction to $(A_1,\deopo)$.

\begin{proof}
Let $(X,\te_1)$ be a covariant representation for the coaction
$\eta_1$ on the \cst-$B$-module
$\E$, consisting of a corepresentation $X \in \cL(A_1 \ot \E)$ and a
representation $\te_1 : C \recht \cL(\E)$ satisfying the relation
$$(\io \ot \te_1) \etatil(y) = X (1 \ot \te_1(y)) X^* \quad\text{for
  all}\;\; y \in C \; .$$
We induce the corepresentation $X$ to a corepresentation $Y := \Ind X$
  on the induced \cst-$B$-module $\cF := \Ind \E$. Recall that
\begin{equation} \label{eq.recall}
H \ot \cF \cong \cI \rot{\pil} (H \ot \E) \; ,
\end{equation}
where the representation $\pil : \Mh_1 \recht \cL(H \ot \E)$ is
determined by $(\io \ot \pil)(W_1) = (\io \ot \pih)(W_1)_{12}
X_{13}$. Using the covariance of $\te_1$, it is easily checked that
$\pil(\Mh_1)$ and $(\io \ot \te_1)\etatil(C)$ commute. So, we get a
well defined representation
$$\tetil : C \recht \cL(H \ot \F) : \tetil(y) (v \rot{\pil} \xi) =
v \rot{\pil} (\io \ot \te_1)\etatil(y) \xi \quad\text{for all}\; \; v
\in  \cI, \xi \in H \ot \E \; .$$
We implicitly used the identification in \eqref{eq.recall}
when defining $\tetil$.
Since $(\io \ot \te_1)\etatil(C)$ commutes with $M' \ot 1$, it follows
that $\tetil(C)$ commutes with $M' \ot 1$ as well. By definition of
$\tetil$, we have
$$(\io \ot \tetil)\etatil(y) = V_{21} (1 \ot \tetil(y)) V_{21}^*
\quad\text{for all}\;\; y \in C \; .$$
Finally, we have $\pil : \Mh \ltimes Q \recht \cL(H \ot \F)$ and the
images of $\pil(\Mh \ltimes Q)$ and $\tetil(C)$ commute. This means
in particular that $1 \ot \tetil(C)$ commutes with $W_{12} Y_{13}$.

It follows that, for all $y \in C$,
$$W^*_{12} (Y^* \tetil(y) Y)_{23} W_{12} = Y^*_{23} Y^*_{13} W^*_{12} (\io
\ot \tetil(y)) W_{12} Y_{13} Y_{23} = (Y^* \tetil(y) Y)_{23} \; .$$ Hence, $Y^* \tetil(y) Y$ commutes with $\Mh \ot 1$. We already know that $Y^*
\tetil(y) Y$ commutes with $M' \ot 1$. So, we find a representation $\te : C \recht \cL(\F)$ such that $\tetil(y) = Y (1 \ot \te(y)) Y^*$. From the
relations stated above, we conclude that the image of $\te$ commutes with the image of the homogeneous space $\rho : D \recht \cL(\F)$ and that $(\io
\ot \te)\etatil(y) = Y (1 \ot \te(y))Y^*$ for all $y \in C$.

So, we have found a triple $(\rho,Y,\te)$, which means a representation of $D \ltimesfull \Ahu \ltimesfull C$ on $\F = \Ind \E$.

We can perform as well the inverse induction. Let $\F$ be a
\cst-$B$-module. Suppose that $(\rho,Y,\te)$
is a triple consisting of commuting representations $\rho$ of $D$ and $\te$ of
$C$ on $\F$, and a corepresentation $Y \in \cL(A \ot \F)$ satisfying the
covariance relations in \eqref{eq.againcov}.

As in the proof of Theorem \ref{thm.firstimprim}, we find the
\cst-$B$-module $\E$ by $H \ot \E \cong \cI^* \rot{\pil} (H \ot \F)$,
where $\pil : \Mh \ltimes Q \recht \cL(H \ot \F)$ is the strict
$^*$-homomorphism determined by
$$(\io \ot \pil)(W) = W_{12} Y_{13} \quad\text{and}\quad \pil(x) = 1
\ot \rho(x) \quad\text{for all} \;\; x \in D \; .$$
We find a corepresentation $X \in \cL(A_1 \ot \E)$ such that $Y= \Ind
X$.

Defining $\tetil : C \recht \cL(H \ot \F) : \tetil(y) = Y(1 \ot
\te(y))Y^*$, it is clear that the images of $\tetil$ and $\pil$
commute. So, we obtain a representation $\tetil_1 : C \recht \cL(H \ot
\E)$ given by
$$\tetil_1(y) (v \rot{\pil} \xi) = v \rot{\pil} \tetil(y)\xi
\quad\text{for all}\;\; v \in \cI^*, \xi \in H \ot \F \; .$$
Exactly as it was the case in the induction procedure above, we find
that $\tetil_1(C)$ commutes with $M' \ot 1$ and satisfies
\begin{equation}\label{eq.weerwat}
V_{21} (1 \ot \tetil_1(y)) V^*_{21} = (\io \ot \tetil_1)\etatil(y)
\quad\text{for all}\;\; y \in C \; .
\end{equation}
Finally, $1 \ot \tetil_1(C)$ commutes with $(\io \ot \pih)(W_1)_{12}
X_{13}$.

Suppose now first that $\eta$ is a maximal coaction. This means that the natural surjective $^*$-homomorphism $\Au \ltimesfull \Ahu \ltimesfull C
\recht \K(H) \ot C$ is an isomorphism. At the end of the previous paragraph, we found a representation $\tetil_1$ of $C$ on $H \ot \E$, which
commutes with $M' \ot 1$ and satisfies the covariance relation \eqref{eq.weerwat} with respect to $\Mh' \ot 1$. The maximality of $\eta$ implies that
there exists a unique representation $\te_1 : C \recht \cL(\E)$ such that $\tetil_1 = (\io \ot \te_1)\etatil$. It follows easily that $\te_1$ is
covariant with respect to $X$. So, we have found a covariant pair $(X,\te_1)$ such that $(\rho,Y,\te)$ is the induction of $(X,\te_1)$. In the same
way as we have shown Theorem \ref{thm.morita}, it follows that there exists a canonical Morita equivalence
$$D \ltimesfull \Ahu \ltimesfull C \Morita \Ahu_1 \ltimesfull C \; .$$
More concretely, the Morita equivalence can be written as $\Ind (\Ahu_1 \ltimesfull C)$.

Suppose next that $\eta$ is an arbitrary continuous coaction.
Consider $\Ah_1 \ltimesred C$ as a \cst-$\Ah_1 \ltimesred C$-module
and define $\F = \Ind(\Ah_1 \ltimesred C)$. Exactly as before the
proof of Theorem \ref{thm.inductionrestriction} we realize $\F$
concretely as
$$\F \cong [(\cI_0 \ot 1) \eta_1(C)] \subset \cL(H_1 \ot C,H \ot C) \;
.$$
As in the proof of Theorem \ref{thm.inductionrestriction}, but now
using the continuous coaction $\cI_0 \recht \M(\cI_0 \ot \Ah)$ given
by $v \mapsto \Vh (v \ot 1) (\io \ot \pih)(\Vh_1^*)$, we find that
$[(\cI_0 \ot 1) \eta_1(C)] = [\eta(C) (\cI_0 \ot 1)]$. It follows that
$$\K(\F) \cong [(\cI_0 \ot 1) \eta_1(C) (\cI_0^* \ot 1) \eta(C)] =
[(\cI_0 \cI_0^* \ot 1) \eta(C)] = [(D \Ah \ot 1) \eta(C)] \; .$$
So we found a canonical Morita equivalence $\Ah_1 \ltimesred C \Morita
[(D \Ah \ot 1) \eta(C)]$. Suppose that $\eta$ is reduced. Then
$$\Vh_{12} \; (\io \ot \eta)[(D \Ah \ot 1)\eta(C)] \; \Vh_{12}^* = [(D \ot 1 \ot 1)
(\deh(\Ah) \ot 1) (1 \ot \eta(C))] = D \ltimesred \Ah \ltimesred C$$
and we are done.
\end{proof}

\begin{remark}
Observe that we have shown that for arbitrary continuous coactions
$\eta : C \recht \M(JAJ \ot C)$, admitting a lift to the universal
level, there is a canonical Morita equivalence
$$\Ah_1 \ltimesred C \Morita [(D\Ah \ot 1)\eta(C)] \; .$$
\end{remark}

\section{Final remarks}

The particular case of inducing a unitary corepresentation of a closed quantum subgroup has been treated by Kustermans \cite{JK2}. His approach, in
the spirit of Mackey, does not allow to prove \cst-algebraic imprimitivity theorems. Of course, one can verify that his induction is unitarily
equivalent to ours. Nevertheless one needs to use the complete machinery of modular theory to prove this result. This is not surprising: the induced
corepresentation of Kustermans involves the canonical implementation (in the sense of \cite{V}) of the coaction of $(M,\de)$ on $Q$ and this is
essentially an object in modular theory. The key result that one has to prove is that the induction of the trivial corepresentation of $(M_1,\de_1)$
is exactly this unitary implementation.

From the naturality and functoriality of our induction procedure, it follows immediately that there is a theorem on \emph{induction in stages}: if
$(M_2,\de_2)$ is a closed quantum subgroup of $(M_1,\de_1)$ and the latter is a closed quantum subgroup of $(M,\de)$, then inducing first from
$(M_2,\de_2)$ to $(M_1,\de_1)$ and then from $(M_1,\de_1)$ to $(M,\de)$ is the same as inducing from $(M_2,\de_2)$ to $(M,\de)$.

In the case where the closed quantum subgroup $(M_1,\de_1)$ of $(M,\de)$ is \emph{normal}, i.e.\ when we have a short exact sequence $e
\longrightarrow (M_2,\de_2) \longrightarrow (M,\de) \longrightarrow (M_1,\de_1) \longrightarrow e$, it follows immediately from the uniqueness
statement in Theorem \ref{thm.existqhs} that the quantum homogeneous space is exactly the reduced \cst-algebra of the quantum group $(M_2,\de_2)$.

The latter example shows moreover that a quantum homogeneous space satisfying the conditions in Theorem \ref{thm.existqhs} may exist even in the
non-regular or non-semi-regular case.

\section{Appendix: \cst- and \wst-modules and their coactions} \label{sec.appendix}

\subsection{Coactions on \cst-modules}

We briefly recall from \cite{KKS} the notion of a coaction on a
Hilbert \cst-module.

\begin{notation}
Let $\E$ be a \cst-$B$-module. Then we denote
$$\M(\E) = \cL(B,\E) \; .$$
\end{notation}

\begin{definition}
Let $\al_B : B \recht \M(B \ot A)$ be a coaction of $(A,\de)$ on $B$
and let $\E$ be a \cst-$B$-module. A coaction of $(A,\de)$ on $\E$
compatible with $\al_B$ is a linear map
$$\al_\E : \E \recht \M(\E \ot A)$$
satisfying
\begin{enumerate}
\item $\al_\E(vx) = \al_\E(v) \; \al_B(x) \quad\text{for all}\;\; v
  \in \E, x \in B \; ;$ \\ $\la \al_\E(v) , \al_E(w) \ra = \al_B(\la
  v,w \ra) \quad\text{for all}\;\; v,w \in \E \; ; $
\item the linear span of $\al_\E(\E)(B \ot A)$ is dense in $B \ot A$;
\item $(\al_\E \ot \io)\al_\E = (\io \ot \de)\al_\E$ (which makes
  sense because of 1) and 2), see \cite{KKS} for details).
\end{enumerate}
\end{definition}

Let $\al_\E$ be a coaction of $(A,\de)$ on $\E$ compatible with
$\al_B$. Then we construct a unitary operator
$$\cV : \E \rot{\al_B} (B \ot A) \recht \E \ot A : \cV (v \rot{\al_B}
x) = \al_\E(v)x \; .$$
It is easy to verify that $\cV$ satisfies the relation
\begin{equation}\label{eq.Vcoaction}
(\cV \rot{\C} \io)(\cV \rot{\al_B \ot \io} 1) = \cV \rot{\io \ot \de}
1 \; .
\end{equation}
This equality holds in $\cL(\E \rot{(\io \ot \de)\al_B} (B \ot A \ot
A), \E \ot A \ot A)$ and its correct interpretation uses the following
identifications.

\begin{equation}\label{eq.precisemeaning}
\begin{CD}
\E \rot{(\al_B \ot \io)\al_B} (B \ot A \ot A) @= \E \rot{(\io \ot
  \de)\al_B} (B \ot A \ot A) \\
@V{\simeq}VV @VV{\simeq}V \\
\bigl(\E \rot{\al_B} (B \ot A) \bigr) \rot{\al_B \ot \io} (B \ot A \ot
A) @. \bigl(\E \rot{\al_B} (B \ot A) \bigr) \rot{\io \ot \de} (B
\ot A \ot A)  \\
@V{\cV \rot{\al_B \ot \io} 1}VV @VVV \\
(\E \ot A) \rot{\al_B \ot \io} (B \ot A \ot A) @. \mbox{}\hspace{13mm}\downarrow \;\; \cV
\rot{\io \ot \de} 1 \\
@V{\simeq}VV  @VVV \\
\bigl(\E \rot{\al_B} (B \ot A)\bigr) \ot A @. (\E \ot A) \rot{\io
  \ot \de} (B \ot A \ot A) \\
@V{\cV \rot{\C} \io}VV @VV{\simeq}V \\
\E \ot A \ot A @= \E \ot A \ot A
\end{CD}
\end{equation}

Whenever we write the symbol $\simeq$ in this diagram, we mean that
there is a natural identification, not involving the coaction on $\E$.

\begin{proposition}
Let $\al_B : B \recht \M(B \ot A)$ be a coaction and $\E$ a
\cst-$B$-module. Let $\al_\E : \E \recht \M(\E \ot A)$ be a linear map.
Then, the following conditions are equivalent.
\begin{enumerate}
\item $\al_\E$ is a coaction of $(A,\de)$ on $\E$ compatible with
  $\al_B$.
\item There exists a coaction on the link algebra $\K(\E \oplus B)$
  that coincides with $\al_\E$ on $\E$ and $\al_B$ on $B$.
\item The formula $\cV(v \rot{\al_B} x) = \al_\E(v)x$ defines a
  unitary operator $\cV$ in $\cL(\E \rot{\al_B} (B \ot A), \E \ot A)$
  satisfying \eqref{eq.Vcoaction}.
\end{enumerate}
\end{proposition}

\begin{remark} \label{rem.corepcoaction}
Let $\E$ be a \cst-$B$-module. If $\cV$ defines a coaction of
$(A,\de)$ on $\E$ which is compatible with the trivial coaction on
$B$, then $\cV \in \cL(\E \ot A,\E \ot A) = \M(\K(\E) \ot A)$ and as
such, $\cV$ is a corepresentation of $(A,\de)$ in $\K(\E)$.
\end{remark}

\begin{remark}
If $(A,\de)$ is a regular l.c.\ quantum group and $\al_B$ a continuous
coaction, then the associated coaction on the link algebra is
automatically continuous.

Indeed, from the compatibility of $\al_\E$ and $\al_B$, as well as the
continuity of $\al_B$, it follows that $\E = [(\io \ot
\om)\al_\E(\E)]$. One can repeat to prove of Proposition 5.8 in
\cite{BSV} to obtain that $[(1 \ot A) \al_\E(\E)] = \E \ot A$ and then
we are done.
\end{remark}

\subsection{Coactions on \wst-modules}

Usually, a von Neumann algebra $M$ is defined as a \cst-algebra acting
non-degenerately on
a Hilbert space such that one of the following equivalent conditions
holds true: $M$ is weakly closed, the unit ball of $M$ is
strongly-$^*$ closed, $M$ is equal to its bicommutant $M\dpr$.

We define in the same way the notion of a \wst-$M$-module. The reader
should convince himself that the proofs of the following two
propositions are elementary.

\begin{proposition}
Let $M \subset \B(H)$ be a von Neumann algebra and let $\E$ be a
\cst-$M$-module. Then, the following conditions are equivalent.
\begin{enumerate}
\item $\E \rot{M} 1 \subset \B(H,\E \rot{M} H)$ is weakly closed.
\item The unit ball of $\E \rot{M} 1$ is strongly-$^*$ closed.
\item $\E \rot{M} 1 = \{ T \in \B(H,\E \rot{M} H) \mid T x = (1
  \rot{M} x) T \quad\text{for all}\; \; x \in M' \}$.
\item The link algebra $\displaystyle \cL(\E \oplus M) = \begin{pmatrix} \cL(\E) \rot{M} 1 & \E \rot{M} 1
    \\ (\E \rot{M} 1)^* & M \end{pmatrix} \subset \B\bigl( (\E \rot{M}
    H) \oplus H \bigr)$ is a von Neumann algebra.
\end{enumerate}
If one of these conditions holds true, we call $\E$ a \wst-$M$-module.
\end{proposition}

In the following way, we extend the notion of a normal, unital
$^*$-homomorphism to \wst-modules.

\begin{proposition}
Let $M \subset \B(H_M)$ and $N \subset \B(H_N)$ be von Neumann
algebras and $\pi_M : M \recht N$ a normal, unital
$^*$-homomorphism. Let $\E$ be a \wst-$M$-module and $\F$ a \wst-$N$-module.

Suppose that $\pi_\E : \E \recht \F$ is a linear map such that
\begin{itemize}
\item $\pi_\E(vx) = \pi_\E(v) \pi_M(x) \quad\text{for all}\;\; v \in
  \E, x \in M \; ,$
\item $\la \pi_\E(v) , \pi_\E(w) \ra = \pi_M( \la v,w \ra )
  \quad\text{for all}\;\; v,w \in \E \; .$
\end{itemize}

Then, $\pi_\E$ is automatically strongly$^*$ continuous on the unit
ball of $\E$. Moreover, the following conditions are equivalent.
\begin{enumerate}
\item $\pi_\E(\E) H_N := (\pi_\E(\E) \rot{N} 1) H_N$ is dense in $\F
  \rot{N} H_N$.
\item $\pi_\E$ and $\pi_M$ extend to a unital, normal,
  $^*$-homomorphism $\cL(\E \oplus M) \recht \cL(\F \oplus N)$.
\end{enumerate}

If one of these conditions holds true, we say that $\pi_\E$ is a
non-degenerate morphism compatible with $\pi_M$. In that case, the
extension to the link algebra $\cL(\E \oplus M)$ is unique.
\end{proposition}

Having spelled out the notion of a non-degenerate morphism, we can
study coactions on \wst-modules.

Remark that is obvious how to define outer and interior tensor
products of \wst-modules, in the same spirit as for \cst-modules.

\begin{definition}
Let $\al_N : N \recht N \ot M$ be a coaction of a l.c.\ quantum group
$(M,\de)$ on the von Neumann algebra $N$. Let $\E$ be a
\wst-$N$-module.

Let $\al_\E : \E \recht \E \ot M$ be a non-degenerate morphism compatible with
$\al_N$. Then, the following two conditions are equivalent.
\begin{enumerate}
\item $(\al_\E \ot \io)\al_\E = (\io \ot \de)\al_\E \; .$
\item The extension of $\al_\E$ and $\al_N$ to the link algebra
  $\cL(\E \oplus N)$ is a coaction of $(M,\de)$ on the von Neumann
  algebra $\cL(\E \oplus N)$.
\end{enumerate}
In that case, we say that $\al_\E$ is a coaction of $(M,\de)$ on $\E$ compatible with $\al_N$.
\end{definition}

\subsection{The interior tensor product of a \wst- and a \cst-module}

The subtle point of this paper is the construction of an interior tensor product of
a \wst-module and a \cst-module, as well as an interior tensor product
of compatible coactions. We have seen in Definition
\ref{def.stricthomo} the notion of a strict homomorphism $N \recht
\cL(\E)$, where $N$ is a von Neumann algebra and $\E$ a \cst-module.

The two natural constructions to develop next are the following. They
are the crucial technical ingredients for the approach to induction
presented in this paper.

\begin{enumerate}
\item Define the interior tensor product $\I \rot{\pi} \E$
of the \wst-$N$-module $\I$ and the \cst-$B$-module $\E$ when $\pi : N
\recht \cL(\E)$ is strict.
\item Given a coaction of $(M,\de)$ on $\I$ and a coaction of its
  \cst-algebraic version $(A,\de)$ on $\E$ such that $\pi$ is
  covariant, construct an interior tensor product coaction of
  $(A,\de)$ on $\I
  \rot{\pi} \E$ in the spirit of Baaj \& Skandalis \cite{KKS}, who
  deal with the interior tensor product of coactions on \cst-modules.
\end{enumerate}

\begin{definition} \label{def.tensorwst}
Let $\I$ be a \wst-$N$-module and $\E$ a \cst-$B$-module. Let $\pi :
N \recht \cL(\E)$ be a strict $^*$-homomorphism.

Then, the algebraic tensor product $\I \rotalg{\pi} \E$ can be completed
to a \cst-$B$-module $\I \rot{\pi} \E$ using the inner product
$$\la v \rot{\pi} w , v' \rot{\pi} w' \ra = \la w, \pi(\la v,v' \ra)
w' \ra \quad\text{for}\;\; v,v' \in \I, w,w' \in \E \; .$$
\end{definition}

\begin{remark}
Observe that when $(v_i)$ is a bounded net in $\I$ converging
strongly$^*$ to $v \in \I$ and $w \in \E$, then $(v_i \rot{\pi} w)$
converges (in norm) to $v \rot{\pi} w$.
In particular, if $\I_0$ is a subspace of $\I$ whose unit ball is
strongly$^*$ dense in $\I$, then the algebraic tensor product $\I_0
\rotalg{\pi} \E$ is dense in $\I \rot{\pi} \E$.

In the situation of Definition \ref{def.tensorwst}, $\cL(\I)$ is a von
Neumann algebra and the $^*$-homomorphism $\cL(I) \recht \cL(\I
\rot{\pi} \E) : x \mapsto x \rot{\pi} 1$ is a strict $^*$-homomorphism.
\end{remark}

Suppose now that a l.c.\ quantum group $(M,\de)$ is coacting on $N$ by
$\al_N : N \recht N \ot M$. Suppose
that its \cst-algebraic companion $(A,\de)$ is coacting on a \cst-$B$-module
$\E$. In particular, we have the coaction
$$\be_{\cL(\E)} : \cL(\E) \recht \cL(\E \ot A) \; .$$
Suppose that $\pi : N \recht \cL(E)$ is a strict $^*$-homomorphism. We
then want to give a meaning to the covariance relation
$$(\pi \ot \io) \al_N = \be_{\cL(\E)} \pi \; .$$
Using the representation of $A$ on the Hilbert space $H$, we have
$\cL(\E \ot A) \hookrightarrow \cL(\E \ot H)$. Lemma
\ref{lem.extendstrict} below
tells how to define a strict $^*$-homomorphism $\pi \ot \io : M \ot
\B(H) \recht \cL(\E \ot H)$. This leads to the following definition.

\begin{definition} \label{def.covariant}
In the situation described in the previous paragraph, we say that
$\pi$ is covariant when the equation $(\pi \ot \io)\al_N(x) =
\be_{\cL(\E)} \pi(x)$ holds in $\cL(\E \ot H)$ for all $x \in N$.
\end{definition}

\begin{lemma} \label{lem.extendstrict}
Let $N$ be a von Neumann algebra and $\E$ a \cst-$B$-module. Suppose
that $\pi : N \recht \cL(E)$ is a strict $^*$-homomorphism. Let $H$ be
a Hilbert space. Then,
there exists a unique strict $^*$-homomorphism $\pi \ot \io : N \ot
\B(H) \recht \cL(\E \ot H)$ satisfying $(\pi \ot \io)(x \ot y) =
\pi(x) \ot y$.
\end{lemma}
\begin{proof}
Consider the \wst-$N$-module $N \ot H$. Identifying $\E \ot H \simeq
(N \ot H) \rot{\pi} \E$, we get a strict $^*$-homomorphism $\cL(N \ot
H) \recht \cL(\E \ot H)$. It is clear that $\cL(N \ot H) = N \ot
\B(H)$ and so, we are done.
\end{proof}

We finally want to construct the interior tensor product of a coaction
on a \wst-module and a coaction on a \cst-module, following Baaj \&
Skandalis \cite{KKS} who made the interior tensor product of coactions
on \cst-modules.

We fix the following data.
\begin{itemize}
\item
Let $\al_N : N \recht N \ot M$
be a coaction of a l.c.\ quantum group $(M,\de)$ on a von Neumann
algebra $N$. Let $\I$ be a \wst-$N$-module and $\al_\I : \I \recht \I
\ot M$ a compatible coaction of $(M,\de)$ on $N$.
\item
Let $(A,\de)$ be the \cst-algebraic companion of $(M,\de)$.
Let $\be_B : B \recht \M(B \ot A)$ be a coaction of $(A,\de)$ on the
\cst-algebra $B$. Let $\E$ be a \cst-$B$-module equipped with a
compatible coaction $\be_\E : \E \recht \M(\E \ot A)$.
\item
Let $\pi : N \recht \cL(\E)$ be a strict $^*$-homomorphism which is
covariant in the sense of Definition \ref{def.covariant}.
\end{itemize}

\begin{proposition} \label{prop.productofcoactions}
In the situation above, there exists a unique coaction $\ga_\F$ of $(A,\de)$ on
$\F:=\I \rot{\pi} \E$ compatible with $\be_B$ and satisfying
\begin{equation}\label{eq.productcoaction}
\ga_\F(v \rot{\pi} w) x = \al_\I(v) \rot{\pi \ot \io}
\bigl(\be_\E(w) x \bigr) \quad\text{for all}\;\; v \in \I, w \in \E, x
\in B \ot \K(H) \; .
\end{equation}
Here we use the embeddings $\M(\F \ot A) \hookrightarrow \M(\F \ot
\K(H))$ and $\M(\E \ot A) \hookrightarrow \M(\E \ot \K(H))$ to give a
meaning to the previous equality.
\end{proposition}

Remark that is included in the contents of the proposition that
$\ga_\F : \F \recht \M(\F \ot A)$, which is non-obvious from the
defining relation \eqref{eq.productcoaction}.

\begin{proof}
We shall write $\K$ for $\K(H)$ throughout the proof.
We define a unitary $$\cV \in \cL((\I \rot{\pi} \E) \rot{\be_B} (B \ot
\K), (\I \rot{\pi} \E) \ot \K)$$ by the formula
$$\cV\bigl( (v \rot{\pi} w) \rot{\be_B} x \bigr) = \al_\I(v) \rot{\pi
  \ot \io} \bigl(\be_\E(w) x\bigr) \; .$$
The slightly non-trivial point to check is the surjectivity of
  $\cV$. To prove this, it suffices to check that any element of $\I
  \ot \B(H)$ can be approximated in the strong$^*$ topology by a
  bounded net in $\operatorname{span}\bigl(\al_\I(\I)(1 \ot
  \K)\bigr)$. This last result follows from the fact that $\al_I$ and
  $\al_N$ combine to a coaction of $(M,\de)$ on the link algebra, on
  which we can apply the results of \cite{V}.

We claim that $\cV$ satisfies the relation
\begin{equation}\label{eq.wat}
(\cV \rot{\C} \io)(\cV \rot{\be_B \ot \io} 1) = \cV \rot{\io \ot
  \de} 1
\end{equation}
which holds in $\cL((\I \rot{\pi} \E) \rot{(\be_B \ot \io)\be_B} (B
  \ot \K \ot \K), (\I \rot{\pi} \E) \ot \K \ot \K)$ and which should
  be given a precise meaning as follows. Using the multiplicative
  unitary $W \in \M(A \ot \hat{A})$, we define $\de : \K \recht \M(\K
  \ot \K) : \de(k) = W^* (1 \ot k) W$, extending the comultiplication
  $\de$ on $A$. Then Equation \eqref{eq.wat} gets a precise meaning as
  in Equation \eqref{eq.precisemeaning} replacing systematically $A$
  by $\K$.

Let $(e_i)$ be an approximate unit in $B \ot \K$. Then, for $v \in
\I$, $w \in \E$ and $x \in B \ot \K \ot \K$, we have
$$(\cV \rot{\be_B \ot \io} 1) \bigl( (v \rot{\pi} w) \rot{(\be_B \ot
  \io)\be_B} x \bigr) = \lim_i (\cV \rot{\be_B \ot \io} 1)\bigl(
  \bigl( (v \rot{\pi} w) \rot{\be_B} e_i \bigr) \rot{\be_B \ot \io} x
  \bigr)
= \lim_i \bigl( \al_I(v) \rot{\pi \ot \io} \be_\E(w)e_i \bigr)
  \rot{\be_B \ot \io} x \; .$$
We now consider
\begin{align*}
& (\I \ot M) \rot{\pi \ot \io} (\E \ot \K) \rot{\be_B \ot \io} (B \ot
\K \ot \K) \simeq \bigl( (\I \rot{\pi} \E) \rot{\be_B} (B \ot
\K)\bigr) \ot \K \\ \overset{\cV \rot{\C} \io}{\recht} & (\I \rot{\pi} \E)
\ot \K \ot \K \simeq (\I \ot M \ot M) \rot{\pi \ot \io \ot \io} (B \ot
\K \ot \K) \; .
\end{align*}
This chain of maps applied to an elementary tensor yields
\begin{align*}
& (v \ot a) \rot{\pi \ot \io} (w \ot k) \rot{\be_B \ot \io} (y \ot l)
  \\
\mapsto & \bigl( (v \rot{\pi} w) \rot{\be_B} y \bigr) \ot akl \\
\mapsto & (\al_\I(v) \rot{\pi \ot \io} \be_\E(w)y) \ot akl \\
\mapsto & (\al_\I \ot \io)(v \ot a) \rot{\pi \ot \io \ot \io} \bigl(
  (\be_\E \ot \io)(w \ot k) (y \ot l) \bigr) \; .
\end{align*}

It follows that
\begin{align*}
(\cV \rot{\C} \io)(\cV \rot{\be_B \ot \io} 1) \bigl( (v \rot{\pi} w) \rot{(\be_B \ot
  \io)\be_B} x \bigr) &= \lim_i (\al_I \ot \io)\al_I(v) \rot{\pi \ot
  \io \ot \io} (\be_\E \ot \io)(\be_\E(w)e_i) x \\ &= (\io \ot
  \de)\al_\I(v) \rot{\pi \ot \io \ot \io} (\io \ot \de)\be_\E(w) x \;
  .
\end{align*}

On the other hand, we have
$$(\cV \rot{\io \ot \de} 1) \bigl( (v \rot{\pi} w) \rot{(\be_B \ot
  \io)\be_B} x \bigr) = \lim_i (\cV \rot{\io \ot \de} 1)\bigl(
  \bigl( (v \rot{\pi} w) \rot{\be_B} e_i \bigr) \rot{\io \ot \de} x
  \bigr)
= \lim_i (\al_\I(v) \rot{\pi \ot \io} \be_\E(w)e_i) \rot{\io \ot \de}
  x \; .$$
We identify
$$\bigl( (\I \ot M) \rot{\pi \ot \io} (\E \ot \K) \bigr) \rot{\io \ot
  \de} (B \ot \K \ot \K) \simeq (\I \ot M \ot M) \rot{\pi \ot \io \ot
  \io} (B \ot \K \ot \K)$$
which is given by
$$(v \rot{\pi \ot \io} w) \rot{\io \ot \de} x \mapsto (\io \ot \de)(v)
\rot{\pi \ot \io \ot \io} (\io \ot \de)(w)x \; .$$
Hence, we conclude that
$$(\cV \rot{\io \ot \de} 1) \bigl( (v \rot{\pi} w) \rot{(\be_B \ot
  \io)\be_B} x \bigr) = \lim_i (\io \ot \de)\al_\I(v) \rot{\pi \ot \io
  \ot \io} (\io \ot \de)(\be_\E(w)e_i) x = (\io \ot \de)\al_\I(v) \rot{\pi \ot \io
  \ot \io} (\io \ot \de)(\be_\E(w)) x\; .$$
This proves our claim. From the lemma following this proposition, we
  get that $\cV \in \cL((\I \rot{\pi} \E) \rot{\be_B} A, (\I \rot{\pi}
  \E) \ot A)$ and that we get the desired coaction $\ga_\F$ on $\F =
  \I \rot{\pi} \E$ as stated in the proposition.
\end{proof}

\begin{lemma}
Let $\E$ be a \cst-$B$-module and let $\be_B : B \recht \M(B \ot A)$
be a coaction of a l.c.\ quantum
group $(A,\de)$ on $B$. Suppose that $\cV \in \cL\bigl(\E \rot{\be_B} (B
\ot \K(H)), \E \ot \K(H)\bigr)$ is a unitary satisfying
$$(\cV \rot{\C} \io)(\cV \rot{\be_B \ot \io} 1) = \cV \rot{\io \ot
  \de} 1 \; .$$
Here, we denote $\de : \K(H) \recht \M(\K(H) \ot \K(H)) : \de(k) = W^*
  (1 \ot k) W$ and we refer to \eqref{eq.precisemeaning} for
  the precise meaning of the formula satisfied by $\cV$.

Then, $\cV \in \cL\bigl(\E \rot{\be_B} (B
\ot A), \E \ot A\bigr)$ and hence, there exists a unique coaction
  $\be_\E$ of $(A,\de)$ on $\E$ compatible with $\be_B$ and satisfying
$$\be_B(v)x = \cV(v \rot{\be_B} x) \quad\text{for all}\;\; v \in \E, x
\in B \ot A \; .$$
\end{lemma}
\begin{proof}
Since $W \in \M(A \ot \K(H))$, we get that in fact $\de : \K \recht
\M(A \ot \K)$. From this, we conclude that
$$\cV \rot{\io \ot \de} 1 \in \cL(\E \rot{(\be_B \ot \io)\be_B} (B \ot
A \ot \K),\E \ot A \ot \K) \; .$$
Further, $\be_B : B \recht \M(B \ot A)$, from which we get that
$$\cV \rot{\be_B \ot \io} 1 \in \cL(\E \rot{(\be_B \ot \io)\be_B} (B
\ot A \ot \K), (\E \rot{\be_B} (B \ot A)) \ot \K) \; .$$
But then, the formula satisfied by $\cV$ guarantees that
$$\cV \rot{\C} \io \in \cL((\E \rot{\be_B} (B \ot A)) \ot \K,\E \ot A
\ot \K) \; .$$
Hence, we find that $\cV \in \cL\bigl(\E \rot{\be_B} (B
\ot A), \E \ot A\bigr)$
\end{proof}

\end{document}